\documentclass[12pt]{amsproc} 
\usepackage{latexsym}

\copyrightinfo{}{}

\newtheorem{theorem}{Theorem \thesection.\!\!} 
\newtheorem*{mtheorem}{Main Theorem} 
\newtheorem{mntheorem}[theorem]{Main Theorem \thesection.\!\!} 
\newtheorem{lemma}[theorem]{Lemma \thesection.\!\!} 
\newtheorem{vanlemma}[theorem]{Vanishing Lemma \thesection.\!\!} 
\newtheorem{proposition}[theorem]{Proposition \thesection.\!\!} 
\newtheorem{corollary}[theorem]{Corollary \thesection.\!\!} 
 
\theoremstyle{definition} 
\newtheorem{definition}[theorem]{Definition \thesection.\!\!} 
 
\newtheorem{remark}[theorem]{Remark \thesection.\!\!} 
\newtheorem{example}[theorem]{Example \thesection.\!\!} 
 
\newtheorem*{prexao}{Proof of Example 1.5} 
\newtheorem*{procoo}{Proof of Corollary 4.5} 
\newtheorem*{protheo}{Proof of Theorem 4.4} 
 
\DeclareMathOperator\Hom{Hom} 
\DeclareMathOperator\coker{coker} 
\DeclareMathOperator\Tor{Tor} 
\DeclareMathOperator\pt{pt} 
\def\A{{\mathcal A}} 
\def\E{{\mathcal E}} 
\def\F{{\mathcal F}} 
\def\IH{{\mathcal IH}} 
\def\H{{\mathcal H}} 
\def\b{^{\scriptscriptstyle\bullet}} 
\def\mal{\mathbin{\!\cdot\!}} 
\def\CC{{\mathbb C}} 
\def\TT{{\mathbb T}} 
\def\QQ{{\mathbb Q}} 
\def\ZZ{{\mathbb Z}} 
\def\PP{{\mathbb P}} 
\def\NN{{\mathbb N}} 
\def\RR{{\mathbb R}} 
\def\mm{{\mathbb m}} 
\def\mm{{\mathfrak m}} 
\def\cld{{\rm cld}} 
\def\longto{\longrightarrow} 
 
\def\mimapdown#1{\big\downarrow\rlap{$\vcenter{\hbox{$\scriptstyle#1$}}$}} 
 
\def\mimapup#1{\big\uparrow\rlap{$\vcenter{\hbox{$\scriptstyle#1$}}$}} 
 
\begin{document} 
\typeout{Topmatter} 
 
\title[Equivariant Intersection Cohomology of Toric Varieties] 
{Equivariant Intersection Cohomology\\ of Toric Varieties} 
 
\author[G. Barthel]{Gottfried Barthel} 
\address{G.B.: Fakult\"at f\"ur Mathematik und 
Informatik, Universit\"at Konstanz, Fach D 203, D-78457 Konstanz} 
\email{Gottfried.Barthel@uni-konstanz.de} 
\author[J.-P. Brasselet]{Jean-Paul Brasselet} 
\address{J.-P.B.: IML - CNRS, Case 907 - Luminy, F-13288 Marseille Cedex~9} 
\email{jpb@iml.univ-mrs.fr} 
\author[K.-H. Fieseler]{\\ Karl-Heinz Fieseler} 
\address{K.-H.F.: Matematiska Institutionen, Box 480, Uppsala 
Universitet, SE-75106 Uppsala} 
\email{khf@math.uu.se} 
\author[L. Kaup]{Ludger Kaup} 
\address{L.K.: Fakult\"at f\"ur Mathematik und Informatik, 
Universit\"at Konstanz, Fach D 203, D-78457 Konstanz} 
\email{Ludger.Kaup@uni-konstanz.de} 
 
 
\dedicatory{Respectfully dedicated to Friedrich Hirzebruch} 
 
\begin{abstract} We investigate the equivariant intersection cohomology 
of a toric variety. Considering the defining fan of the variety as a 
finite topological space with the subfans being the open sets (that 
corresponds to the ``toric'' topology given by the invariant open 
subsets), equivariant intersection cohomology provides a sheaf (of 
graded modules over a sheaf of graded rings) on that ``fan space''. We 
prove that this sheaf is a {\em ``minimal extension sheaf''}, i.e., 
that it satisfies three relatively simple axioms which are known to 
characterize such a sheaf up to isomorphism. In the verification of the 
second of these axioms, a key role is played by ``equivariantly 
formal'' toric varieties, where equivariant and ``usual'' 
(non-equivariant) intersection cohomology determine each other by 
K\"unneth type formul\ae. Minimal extension sheaves can be constructed 
in a purely formal way and thus also exist for non-rational fans. As a 
consequence, we can extend the notion of an {\em equivariantly formal 
fan} even to this general setup. In this way, it will be possible to 
introduce ``virtual'' intersection cohomology for equivariantly formal 
non-rational fans. 
\end{abstract} 
 
\maketitle 
 
 
\tableofcontents 
 
\typeout{Introduction} 
\section*{Introduction} 
 
\noindent The relation between algebraic geometry and topology is 
especially close in the case of smooth compact toric varieties: The 
canonical homomorphism from the Chow-theoretic intersection ring to 
the integral homology intersection ring is an isomorphism. Even when 
we allow quotient singularities, the result remains valid if we replace 
integer with rational coefficients. \par

As is well known\footnote{For the basic theory of toric varieties, 
we refer the reader to the pertinent monographs \cite{Ew}, \cite{Fu}, 
or \cite{Oda$_1$}.}, toric varieties admit a description in terms of 
combinatorial-geometric objects: cones, fans, and polytopes. In the 
case considered above, these combinatorial-geometric objects allow us 
to compute the intersection ring explicitely. In the opposite direction, 
properties of the intersection ring known from algebraic geometry 
yield consequences for the describing geometric object. A most striking 
application is Stanley's beautiful proof -- using the hard Lefschetz 
theorem for projective toric varieties -- for the necessity of the 
conditions that characterize the face numbers of simplicial polytopes, 
conjectured by \hbox{McMullen} (see, e.g., \cite[\S5.6]{Fu}). That 
example should suffice to explain the great interest of toric 
varieties in studying the combinatorial geometry of fans -- as a 
variation of the title of Hirzebruch's landmark book, one might call 
this an {\it Application of Topological Methods from Algebraic Geometry 
to Combinatorial Convexity}. \par
 
Stated in terms of the defining fan, a toric variety $X = X_{\Delta}$ 
is compact and $\QQ$-smooth (i.e., a rational homology manifold) if the 
fan~$\Delta$ is {\it complete\/} and {\it simplicial}. When studying the 
(rational) cohomology algebra of such a toric variety~$X$, a crucial 
role is played by the (rational) {\it Stanley-Reisner ring\/} 
$$ 
  S_\Delta := \QQ[\Delta^1] / I \;, 
$$ 
where $\Delta^1$ is the set of {\it rays\/} (i.e., one-dimensional 
cones) of~$\Delta$, and $\QQ[\Delta^1]$ is the polynomial algebra 
$\QQ\big[(t_\rho)_{\rho \in \Delta^1}\big]$ on free 
generators~$t_{\rho}$ that are in one-to-one correspondence to the 
rays~$\rho$, and where~$I$ is the homogeneous ideal 
$$ 
I := \Bigl\langle \prod_{j=1}^k t_{\rho_{j}} \;;\; \sum_{j=1}^k 
\rho_{j} \not\in \Delta \Bigr\rangle 
$$ 
generated by those square-free monomials where the rays corresponding 
to the factors do not span a cone of the fan. The theorem of Jurkiewicz 
and Danilov \cite[3.3, p.~134]{Oda$_1$} describes the (rational) 
cohomology algebra $H\b(X)$ as the quotient of $S_\Delta$ modulo the 
ideal~$\tilde J$, the image of the homogeneous ideal 
$$ 
J := \Bigl\langle \sum_{\rho \in \Delta^1} \chi(v_\rho) \mal t_\rho 
\;;\; \chi \in \Hom(\TT, \CC^*) \Bigr\rangle \subset \QQ[\Delta^1]\;, 
$$ 
where $v_\rho$ denotes the unique primitive lattice vector spanning 
the ray~$\rho$, and $\TT$, the acting torus. We thus have isomorphisms 
$$ 
H\b(X_{\Delta},\QQ) \,\cong\, S_\Delta/\tilde J \,\cong\,
\QQ[\Delta^1]\, /\, (I + J)\;. \leqno{(\dagger)} 
$$ \par 
 
The ideals~$I$ and~$J$ are most naturally interpreted using the 
one-to-one correspondence between rays and $\TT$-invariant irreducible 
divisors on~$X$: The monomials generating~$I$ correspond to sets of 
divisors with empty intersection, and the generators of~$J$, to 
$\TT$-invariant principal divisors. In fact, identifying a product of 
generators with the intersection of the corresponding divisors, this 
correspondence extends to isomorphisms between the rings 
$\QQ[\Delta^1]/(I+J) \cong H\b(X)$ and the rational Chow ring of~$X$. 
-- We note that the ideal~$I$, and hence the Stanley-Reisner ring, 
depends only on the combinatorial structure of the fan. \par 
 
The isomorphism $\QQ[\Delta^1]/(I+J) \cong H\b(X)$ is a morphism of 
graded rings that multiplies degrees by two. As a consequence, the 
odd-dimensional cohomology of~$X$ vanishes. Furthermore, it turns out 
that the Betti numbers are explicitly computable in terms of 
combinatorial data of the fan. \par 
 
Besides of being a graded algebra that encodes the combinatorics of 
the fan~$\Delta$, the Stanley-Reisner ring~$S_{\Delta}$ has a second, 
more geometric, inter\-pre\-ta\-tion. To give it, we consider, for a 
completely arbitrary fan~$\Delta$, the graded algebra $\A\b(\Delta)$ of 
{\it (rational) $\Delta$-piecewise polynomial functions\/} on the 
support $|\Delta|$ of the fan\footnote{Unless otherwise stated, cones 
and fans are always considered as subsets of the {\it rational\/} 
vector space $N_{\QQ}$ generated by the one parameter subgroups of the 
torus.}, 
i.e., functions $f \colon |\Delta| \to \QQ$ whose restriction 
to each cone $\sigma \in \Delta$ extends to a -- unique -- polynomial 
function $f_{\sigma}$ on the linear span $V_{\sigma} := \sigma + 
(-\sigma)$ of~$\sigma$. Since the algebra of polynomial functions 
on~$V_{\sigma}$ is $S\b(V_{\sigma}^*)$, the symmetric $\QQ$-algebra of 
the dual vector space $V_{\sigma}^*$, we arrive at the formal definition 
$$ 
  \A\b(\Delta) := \{f \colon |\Delta| \to \QQ \;;\; \forall\,{\sigma 
  \in \Delta}\;\; \exists\,{f_{\sigma} \in S\b(V_{\sigma}^*)}: 
  f|_{\sigma} = f_{\sigma}|_{\sigma}\} \;.
$$ 
In addition to this graded $\QQ$-algebra of piecewise polynomial 
functions, we also have to consider the graded $\QQ$-algebra 
$$ 
  A\b = S\b(V^*) \cong \QQ\,[u_{1}, \dots, u_{n}] 
$$ 
of (globally) polynomial functions on $V := M_\QQ$, where $u_{1}, \dots, 
u_{n}$ denotes a basis of $V^* := M_\QQ$ (see section~0 for notations). 
The obvious restriction homomorphism $A\b \to \A\b(\Delta)$, 
$f \mapsto f|_{|\Delta|}$ of graded 
algebras is injective if the fan~$\Delta$ contains at least one 
$n$-dimensional cone, thus making $A\b$ a graded subalgebra of 
$\A\b(\Delta)$. \par 

In the sequel, the algebra $\A\b(\Delta)$ will be used without any 
restriction on the fan. But assuming again that~$\Delta$ is complete and 
simplicial for the moment, we obtain 
a homomorphism $\QQ[\Delta^1] \to \A\b(\Delta)$ of graded $\QQ$-algebras 
by associating to each generator~$t_{\rho}$ of $\QQ[\Delta^1]$ the 
unique piecewise linear function on $|\Delta|$ that takes the value~1 at 
the vector~$v_{\rho}$ and vanishes on all other rays. This homomorphism 
is surjective, and its kernel coincides with the ideal~$I$, so we get 
the identification 
$$ 
  S_\Delta \buildrel \cong \over \longrightarrow 
  \A\b(\Delta)  \leqno{(*)} 
$$ 
that provides the geometric interpretation of the Stanley-Reisner ring.
\par
 
Under that identification, the ideal $\tilde J$ in $S_\Delta$ corresponds 
to the homogeneous ideal $\mm \cdot \A\b(\Delta)$ generated by the 
(globally) linear functions. (In the (sub-)algebra $A\b$, the ideal 
generated by linear functions is the unique homogeneous maximal ideal 
$\mm = A^{>0}$ of polynomial functions vanishing at the origin.) Using 
the above isomorphism~$(\dagger)$, we may look at $H\b(X)$ as quotient 
modulo that ideal. We can thus describe the cohomology algebra of a 
compact, $\QQ$-smooth toric variety $X_{\Delta}$ by an isomorphism 
$$ 
  H\b(X_{\Delta}) 
  \,\cong\, (A\b/\mm) \otimes_{A\b} \A\b(\Delta) \;. 
$$ \par

On the other hand, the ring $\A\b(\Delta)$ itself admits a direct 
topological interpretation: There is a natural isomorphism 
$$ 
H\b_{\TT}(X_\Delta) \buildrel\cong\over\longrightarrow 
\A\b(\Delta)\leqno{(**)} 
$$ 
with the {\it equivariant cohomology algebra\/} of $X_\Delta$ (see, 
e.g., \cite{BriVe}, \cite{GoKoMPh}, and the first section of this 
article). Hence, the two isomorphisms $(*)$ and $(**)$ together yield 
a third, topological, interpretation of the Stanley-Reisner ring. 
\par\medskip
 
The case of {\it non-simplicial\/} fans is much more elusive, as the 
above results do not remain valid. In fact, the cohomology becomes 
quite a delicate object that is difficult to compute. In particular, 
the Betti numbers are not always determined by the combinatorial type 
of the fan.\par

Fortunately, for compact toric varieties with arbitrary singularities, 
intersection cohomology (with respect to ``middle perversity'') is known 
to behave much better than ``usual'' 
cohomology in many respects, e.g., Poincar\'e duality still holds, there 
is a Hard Lefschetz Theorem in the projective case, and intersection 
cohomology Betti numbers are combinatorial invariants. With the crucial 
role of the Stanley-Reisner ring for cohomology and with its various 
interpretations in mind, when looking for a substitute on the 
fan-theoretic side that could play a similar role for intersection 
cohomology, we were lead to investigate the properties of the 
{\it equivariant intersection cohomology\/} $IH\b_\TT(X_\Delta)$. This 
graded vector space is endowed with a natural structure of a graded 
module over the graded ring $H\b_\TT(X_\Delta)$. 
 
The aim of this article is to find some analogue of the combinatorial 
description provided by the above isomorphism~$(**)$ in the case of a 
completely arbitrary fan~$\Delta$, replacing equivariant cohomology with 
equivariant intersection cohomology and keeping possible generalizations 
to the case of non-rational fans in mind. It is convenient to adopt a 
sheaf-theoretic point of view: The finite family consisting of all 
$\TT$-invariant open subsets of the toric variety~$X$ is a topology 
on~$X$, and associating to such a ``toric'' open subset $U \subset X$ the 
graded vector space $IH\b_\TT(U)$ yields a presheaf on this topological 
space. Since open toric subvarieties of $X = X_{\Delta}$ 
are in one-to-one correspondence to subfans of the defining fan~$\Delta$, 
the ``toric'' topology on~$X_{\Delta}$ corresponds to the ``fan topology'' 
on the finite set~$\Delta$, namely, the topology given by the collection 
of the subfans $\Lambda \preceq \Delta$, together with the empty set, as 
family of open sets. For such a subfan~$\Lambda$, we then understand the 
graded $A\b$-module $IH\b_{\TT}(X_\Lambda)$ as the module of sections 
over~$\Lambda$ of the {\it equivariant intersection cohomology 
(pre-)sheaf} $\IH\b_\TT$ on the ``fan space''~$\Delta$. 
 
To be more precise, we show that the presheaf 
$$ 
  \IH\b_\TT \colon \Lambda \mapsto \IH\b_\TT(\Lambda) := 
  IH\b_\TT(X_\Lambda,\QQ)\;, \quad\hbox{for $\Lambda \preceq \Delta$}\;, 
$$ 
is in fact a sheaf of modules over the {\it sheaf of rational piecewise 
polynomial functions on the support of the fan\/} 
$$ 
  \A\b \colon \Lambda \mapsto \A\b (\Lambda) 
  := \{ f: |\Lambda| \longrightarrow \QQ\;;\; 
  \hbox{$\Lambda$-piecewise polynomial}\, \}\; ,
$$ 
and we prove that the sheaf $\IH\b_\TT$ has the three properties stated 
in the following definition. \par 

We call a sheaf $\E\b$ of $\A\b$-modules a {\it minimal extension 
sheaf\/} (of the constant sheaf~$\QQ$) if it satisfies the following 
conditions\footnote{For a sheaf~$\F$ on~$\Delta$ and a cone~$\sigma$, 
we simply write~$\F(\sigma)$ to denote the space of sections
$\F(\langle\sigma\rangle)$ 
over the affine subfan $\langle\sigma\rangle$ generated by~$\sigma$.}: 
\smallskip 
{\def\litem{\par\noindent \hangindent=36pt\ltextindent} 
\def\ltextindent#1{\hbox to \hangindent{#1\hss}\ignorespaces}\it 
\litem{\rm(N)} {\rm Normalization:\/} For the zero cone~$o := \{0\}$, 
  there is an isomorphism $\E\b(o) \cong \A\b(o)\; (= \QQ\b)$. 
\litem{\rm(PF)} {\rm Pointwise Freeness:\/} For each cone $\sigma \in 
  \Delta$, the module~$\E\b(\sigma) := \E\b\bigl(\langle \sigma 
  \rangle\bigr)$ is {\it free} over~$\A\b(\sigma)$. 
\litem{\rm(LME)} {\rm Local Minimal Extension $\bmod\ \mm$:\/} For each 
  cone $\sigma \in \Delta$, the restriction mapping 
$$ 
  \varphi_\sigma \colon \E\b(\sigma) \to \E\b(\partial\sigma) 
$$ 
  induces an isomorphism 
$$ 
  \overline \varphi_\sigma \colon 
  \frac {\E\b(\sigma)}{\mm \mal \E\b(\sigma)}
  \buildrel \cong \over \longrightarrow 
  \frac {\E\b(\partial\sigma)} {\mm \mal \E\b(\partial\sigma)} 
$$ 
  of graded vector spaces.\par} 

\smallskip
Restated in other words, a minimal extension sheaf~$\E\b$ on a 
fan~$\Delta$ is characterized as follows: It is a sheaf 
of graded $\A\b$-modules satisfying the equality $\E\b (o) = \QQ\b$ and 
having the property that for each cone $\sigma \in \Delta$, the 
module~$\E\b(\sigma)$ is free and of minimal rank over~$\A\b(\sigma)$ 
such that the restriction $\varphi_\sigma\colon 
\E\b(\sigma) \to \E\b(\partial\sigma) $ is surjective. 
\par 
 
It is not difficult to see that such a sheaf~$\E\b$ can be constructed 
recursively on the $k$-skeletons $\Delta^{\le k}$ of the fan~$\Delta$ 
in a purely formal way, starting from~$\E^{0}=\QQ$, and that the 
resulting sheaf is {\em uniquely determined up to isomorphism}.\par  
 
Using this notion of minimal extension sheaves, the central result of 
the article can be stated as follows: 
 
\begin{mtheorem} 
The equivariant intersection cohomology sheaf $\IH\b_\TT$ on a fan 
space~$\Delta$ is a minimal extension sheaf of~$\QQ\b$. 
\end{mtheorem} 
 
The fact that $\IH\b_{\TT}$ is a minimal extension sheaf has interesting 
consequences that are discussed in \cite{BBFK$_{2}$}: In particular, it 
leads to a simple proof of the inductive ``Local-Global-Formula'' for 
the computation of intersection cohomology Betti numbers of compact 
toric varieties, but also, {\it mutatis mutandis}, it opens a way to 
introduce the analogous sheaf for arbitrary (not necessarily rational) 
fans in real vector spaces. On the other hand, the ``minimal complexes'' 
of \cite{BeLu} occur now naturally as the ``cellular \v{C}ech complexes'' 
of the sheaf $\IH\b_\TT$.\par

\smallskip 
The article is organized as follows: After fixing our (more or less 
standard) notation, we first recall in section~1 the definition and some 
basic properties of the equivariant cohomology of a toric variety; in 
particular, we look at the {\it presheaf\/} $\H\b_{\TT}$ that it 
defines on the fan space. For later use, we also touch upon the 
equivariant Chern class of toric line bundles. In section~2, we 
introduce the equivariant intersection cohomology and prove that it 
vanishes in odd degrees, which implies that it actually defines a 
{\it sheaf\/} $\IH\b_{\TT}$ on the fan space.\par 
 
The characteristic properties (N), (PF), and (LME) of $\IH\b_{\TT}$, 
formalized in the notion of a minimal extension sheaf, are discussed 
in section~3. Whereas property~(N) is easily seen to hold, the other 
two require considerably more work. Property~(PF) is established in 
section~4 by proving that any {\it contractible affine toric variety 
is equivariantly formal\/} (with respect to intersection cohomology): 
The ``usual'' intersection cohomology $IH\b(X)$ is the quotient of the 
equivariant theory $IH\b_\TT(X)$ modulo the homogeneous $A\b$-submodule 
$\mm \mal IH\b_\TT(X)$. The constructions used in the proof are also 
crucial for proving the remaining property (LME) in section~5. The 
article concludes with a discussion of {\it ``equivariantly formal'' 
fans\/}, i.e., fans~$\Delta$ such that the corresponding toric variety 
$X_{\Delta}$ is equivariantly formal. This property holds both for 
complete fans and for ``affine fans'' consisting of a full dimensional 
cone and its faces; more generally, it holds if and only if the 
non-equivariant intersection cohomology vanishes in odd degrees.
\par \smallskip
 
It is our pleasure to thank Michel Brion, Volker Puppe, and the members 
of the ``Algebraic Group Actions'' seminar at Warsaw University led by 
Andrzej Bia{\l}ynicki-Birula for stimulating discussions and for useful 
comments. We appreciate the pertinent remarks of the referee. 
 
\setcounter{section}{-1} 
 
\typeout{1. Preliminaries} 
\section{Preliminaries} 
\setcounter{theorem}{0} 
 
\noindent We use the following {\it notation\/}: We denote with $\TT 
\cong (\CC^*)^n$ the complex {\it algebraic torus\/} of dimension~$n$, 
with $N := \Hom(\CC^*, \TT) \cong \ZZ^n$ the lattice of its one 
parameter subgroups, and with $M := \Hom(\TT, \CC^*) \cong \ZZ^n$ the 
dual lattice of characters. We recall that there are natural 
isomorphism $\TT \cong N \otimes_{\ZZ} \CC^{*} \cong 
\Hom_{\ZZ}(M, \CC^{*})$. Let~$\Delta$ be a fan in the {\it rational\/} 
vector space $V := N_\QQ := N \otimes_\ZZ \QQ$. For two cones~$\sigma$ 
and~$\tau$ in~$\Delta$, we write $\tau \prec \sigma$ if~$\tau$ is a 
{\it proper\/} face of~$\sigma$, and $\tau \prec_1 \sigma$ if~$\tau$ is 
a {\it facet\/} (i.e., a face of codimension~1) of~$\sigma$. -- The 
symbol $\Lambda \prec \Delta$ indicates that a proper subset $\Lambda$ 
of $\Delta$ actually is a {\it subfan.\/} 
 
A cone $\sigma \in \Delta$ {\it generates\/} the {\it ``affine''} subfan 
$\langle\sigma\rangle := \{\sigma\} \cup \partial\sigma$ consisting of 
$\sigma$ and the {\it boundary subfan\/} $\partial\sigma := \{\tau \in 
\Delta\,;\, \tau \prec \sigma\}$. Occasionally, if there is no risk of 
confusion, we simply write~$\sigma$ instead of~$\langle\sigma\rangle$ 
for notational convenience. \par 
 
For a cone $\sigma \in \Delta$, we denote with $V_\sigma := \sigma + 
(-\sigma)$ the linear span of $\sigma$ in $V$, and with $\TT_\sigma$, 
the subtorus of $\TT$ corresponding to the sublattice $N_\sigma := N 
\cap V_\sigma$ as lattice of its one parameter subgroups (i.e., we have 
$\TT_\sigma \cong N_{\sigma} \otimes_{\ZZ} \CC^{*}$). The associated 
lattice of characters is $M_\sigma := M/N_\sigma^{\perp}$. Let us denote 
with $V_{\sigma}^*$ the dual of $V_{\sigma}$ in $V^*:= M_\QQ := M
\otimes_\ZZ \QQ$. The symmetric $\QQ$-algebra $A\b_\sigma := 
S\b(V_{\sigma}^*)$ is naturally isomorphic to the algebra of polynomial 
functions $f \colon V_\sigma \to \QQ$. The restriction
$f \mapsto f|_\sigma$ thus provides an isomorphism 
$$ 
  A\b_\sigma \buildrel \cong \over \longrightarrow 
  \A\b(\sigma) := \A\b(\langle\sigma\rangle) \leqno{(0.1)} 
$$ 
with the algebra of sections over~$\langle\sigma\rangle$ of the 
sheaf~$\A\b$ of piecewise polynomial functions described in the 
introduction. This is compatible with the induced homomorphisms 
$A\b_\sigma \to A\b_\tau$ and $\A\b(\sigma) \to \A\b(\tau)$ for each 
face $\tau \prec \sigma$. -- In view of the cohomological interpretation 
discussed in section~1, we endow the algebra~$A\b_\sigma$ -- just as we 
did with $A\b := S\b(V^*)$ -- with the grading that doubles the usual 
degrees, so linear polynomials (i.e., elements of $V_\sigma^*$) get the 
degree~$2$, etc. 
 
If $F\b$ is a graded $A\b$-module, we denote with $\overline F\b$ its 
residue class module 
$$ 
  \overline F\b := F\b/(\mm \mal F\b) \,\cong\, \QQ\b 
  \otimes_{A\b} F\b \qquad 
  \hbox{(with $\mm := A^{>0}$ and $\QQ\b := A\b/\mm$),} \leqno{(0.2)} 
$$ 
where $\mm$ is the unique homogeneous maximal ideal of $A\b$, and $\QQ\b$ 
is the graded algebra concentrated in degree zero with~$\QQ^{0} := \QQ$. 
-- By means of the natural 
epimorphism $A\b \to A\b_\sigma$ extending the projection $V^* \to 
V^*_\sigma$ (and corresponding to the restriction of polynomial functions 
from $V:=N_{\QQ}$ to $V_{\sigma}$), every $A\b_\sigma$-module 
$F\b_\sigma$ also is an $A\b$-module, and there is a canonical 
isomorphism $\overline F\b_\sigma \cong F\b_\sigma/(\mm_\sigma \mal 
F\b_\sigma)$.\par 
 
In the sequel, we shall freely use the following basic facts on finitely 
generated graded $A\b$-modules~$F\b$: Given a family $(f_1, \dots, f_r)$ 
of homogeneous elements in $F\b$, it generates $F\b$ over $A\b$ if and 
only if the system of residue classes $(\overline f_1, \dots, \overline 
f_r)$ modulo $\mm \mal F\b$ generates the vector space ${\overline{F}}\b$. 
If $F\b$ is free, then $(f_1, \dots, f_r)$ is part of a basis of $F\b$ 
over $A\b$ if and only if $(\overline f_1, \dots, \overline f_r)$ is 
linearly 
independent over~$\QQ$. Furthermore, every homomorphism $\varphi \colon
F\b \to G\b$ of graded $A\b$-modules induces a homomorphism $\overline 
\varphi \colon {\overline F}\b \to {\overline{G}}\b$ of graded vector 
spaces which is surjective if and only if~$\varphi$ is so. If $F\b$ is 
free, then every homomorphism $\psi \colon {\overline F}\b \to 
{\overline{G}}\b$ can be lifted to a homomorphism $\varphi \colon F\b \to 
G\b$; if both $F\b$ and $G\b$ are free, then $\varphi$ is an isomorphism 
if and only if that holds for~$\overline \varphi$.\par 
 
For the affine toric variety $X_\sigma$ associated to~$\sigma$, 
the torus~$\TT_{\sigma}$ is the {\it isotropy subtorus\/} at any point 
in the unique closed orbit $B_{\sigma}$, and this orbit $B_{\sigma}$ is 
$\TT$-isomorphic to the quotient torus $\TT/\TT_{\sigma}$, looked at as 
a homogeneous space. Any splitting $\TT \cong {\TT\,}' \times \TT_\sigma$ 
of the torus~$\TT$ into the isotropy subtorus $\TT_\sigma$ and a 
complementary subtorus~${\TT\,}' \cong \TT/\TT_{\sigma}$ extends to an 
equivariant splitting of affine toric varieties: The choice of an 
arbitrary base point~$x_{o}$ in the open dense orbit~$B_{o}$ determines 
an embedding of~$\TT$ into~$X_{\sigma}$. Denoting with~$Z_\sigma$ the 
closure of~$\TT_\sigma$ in~$X_\sigma$ with respect to this embedding, we 
obtain an isomorphism 
$$ 
  {\TT\,}' \times Z_\sigma \buildrel\cong\over\longto X_\sigma\,,\; 
  (t,z) \mapsto tz \,. \leqno{(0.3)} 
$$ 
The unique point $x_\sigma \in B_\sigma \cap Z_\sigma$ will sometimes be 
referred to as the {\it distinguished point} in the orbit~$B_\sigma$. We 
remark that~$Z_\sigma$, being the $\TT_{\sigma}$-toric variety 
corresponding to~$\sigma$ considered as $N_{\sigma}$-cone, is 
equivariantly contractible to its fixed point~$x_{\sigma}$, so~$X_\sigma$ 
has the closed orbit~$B_\sigma$ as equivariant deformation retract. -- We 
refer to the above splitting~(0.3) as the {\it affine orbit splitting}. 
\par 
 
\goodbreak 
\section {Equivariant Cohomology of Toric Varieties} 
\setcounter{theorem}{0} 
 
\noindent 
Before proceeding to study the equivariant {\it intersection\/} 
cohomology, we first look at the ``usual'' $\TT$-equivariant cohomology 
$H\b_{\TT}(X)$ of a toric variety $X = X_\Delta$. We recall the 
definition: With respect to a fixed identification $\TT \cong 
(\CC^{*})^n$, a universal $\TT$-bundle is given by the principal 
$(\CC^{*})^n$-bundle 
$$ 
  E\TT :=(\CC^\infty \setminus \{0\})^n \longrightarrow 
  B\TT := (\PP_{\infty})^n \;, 
$$ 
the limit of the finite-dimensional approximations 
$$ 
  E_{m}\TT := (\CC^{m+1} \setminus \{0\})^n 
  \longrightarrow B_{m}\TT := (\PP_m)^n 
$$ 
for $m \to \infty$. One considers the associated bundles 
$$ 
  X_\TT := E\TT \times_{\TT} X \longrightarrow B\TT 
  \qquad\hbox{and}\qquad X_{\TT,m} := 
  E_{m}\TT \times_{\TT} X \longrightarrow B_{m}\TT 
$$ 
with fibre~$X$ and defines the {\it (rational) equivariant cohomology 
algebra\/} of~$X$ as follows: 
$$ 
  H\b_{\TT}(X) := H\b_{\TT}(X, \QQ) := H\b(X_{\TT}, \QQ). 
$$ 
For the homogeneous part of some fixed degree $q \ge 0$, there is the 
description 
$$ 
  H^q_\TT(X, \QQ) \cong \lim_{m \to \infty} H^q(X_{\TT,m}, \QQ) \,.
$$ 
The bundle projection $X_\TT \to B\TT$ makes $H\b_{\TT}(X)$ an algebra 
over the cohomology algebra $H\b(B\TT) \cong \QQ[u_{1}, \dots, u_{n}]$ 
of the classifying space $B\TT = (\PP_{\infty})^n$ of~$\TT$, a polynomial 
algebra on~$n$ free generators of degree~$2$.\par 
 
To relate $H\b_{\TT}(X)$ with the combinatorial data encoded in the 
defining fan~$\Delta$ for~$X$, we first recall the following result. 
 
\begin{lemma} \label{isom1.1} 
For each cone $\sigma$, there are natural isomorphisms 
$$ 
  \A\b(\sigma) \cong A\b_\sigma \cong H\b(B\TT_{\sigma}) \cong 
  H\b_\TT(X_{\sigma}) \;,
$$ 
i.e., they are compatible with maps defined by face relations $\tau 
\prec \sigma$. 
\end{lemma} 
 
\begin{proof} The isomorphism on the left-hand side has been discussed 
in~(0.1). The one in the middle is the special case $V = V_{\sigma}$, 
$\TT = \TT_{\sigma}$ of the isomorphism of graded $\QQ$-algebras 
$$ 
  A\b = S\b(V^*) \buildrel \cong \over \longto H\b(B\TT) 
$$ 
induced from the {\it Chern class homomorphism} 
$$ 
  c\colon M \longrightarrow H^2(B\TT)\,, \quad \chi \mapsto c_1(L_\chi) 
$$ 
that associates to a character $\chi \in M \subset V^*=M_{\QQ}$ the 
first Chern class of the line bundle $L_\chi := E\TT \times_{\TT} 
\CC_\chi \to E\TT \times_\TT \{{\rm pt}\} = B\TT$. Here $\CC_\chi$ 
denotes the one-dimensional $\TT$-module with weight $\chi$. 
 
For a cone $\sigma \in \Delta$ and the corresponding subtorus 
$\TT_\sigma$ of $\TT$, the restriction mapping $f \mapsto f|_{V_\sigma}$ 
and the natural mapping $B\TT_{\sigma} \hookrightarrow B\TT$ induce a 
commutative diagram 
\[ 
\begin{array}{ccc} 
  A\b & \longto & A\b_\sigma \\[3pt] 
  \mimapdown{\cong} & & \mimapdown{\cong}\\[3pt] 
  H\b(B\TT) & \longto & H\b(B\TT_{\sigma}) 
\end{array} 
\] 
since we have $L_{\chi}|_{B\TT_{\sigma}} = L_{\chi|_{\TT_\sigma}}$ 
(the restriction of $L_\chi$ to $B\TT_{\sigma}$ is just the line bundle 
associated to the character $\chi|_{\TT_\sigma}$) and Chern classes 
are functorial, thus proving the naturality of the isomorphism in the 
middle. 
 
We now discuss the isomorphism $H\b(B\TT_{\sigma}) \cong 
H\b_\TT(X_{\sigma})$ on the right hand side: The affine orbit splitting 
$(X_{\sigma},\TT) \cong {\TT\,}' \times (Z_\sigma, \TT_{\sigma})$ of 
(0.3) and the $\TT_{\sigma}$-equivariant contraction $Z_{\sigma} \simeq 
\{x_{\sigma}\}$ onto the distinguished point induce isomorphisms 
\begin{eqnarray*} 
  H\b(B\TT_\sigma) & \cong & 
  H\b(E\TT_\sigma \times_{\TT_\sigma} \{x_{\sigma}\}) \cong
  H\b(E\TT_\sigma 
  \times_{\TT_\sigma} Z_\sigma) \\ 
  & = & H\b_{\TT_{\sigma}}(Z_\sigma) 
  \cong H\b_\TT({\TT\,}' \times Z_\sigma) = H\b_\TT(X_\sigma)\;.
\end{eqnarray*} 
The whole construction is natural with respect to some face relation 
$\tau \prec \sigma$. This is easily seen from any splitting of the torus 
in the form $\TT = \TT' \times \TT_\sigma = {\TT\,}' \times ({\TT\,}'' 
\times \TT_\tau) = ({\TT\,}' \times {\TT\,}'') \times \TT_\tau$, since 
the choices of~${\TT\,}'$ and~${\TT\,}''$ do not play any role. 
\end{proof} 
 
\typeout{1.A Equivariant Cohomology etc.} 
\section*{1.A Equivariant Cohomology as a Presheaf on the Fan} 
 
\noindent On the fan space~$\Delta$ defined in the introduction, we 
now consider the {\it presheaf\/} $\H\b_\TT$ of graded algebras that 
is given by 
$$ 
  \H\b_\TT \colon \Lambda 
  \mapsto H\b_\TT(X_\Lambda) \quad\hbox{for $\Lambda \preceq
  \Delta$}\,. 
$$ 
Inverting the isomorphisms of Lemma~1.1 and using the fact that~$\A\b$ 
clearly is a sheaf on~$\Delta$, we obtain the following result: 
 
\begin{corollary} 
There is a homomorphism of presheaves $\H\b_\TT \to \A\b$ 
that is an isomorphism on the stalks, so $\A\b$ is the associated 
sheaf to the presheaf~$\H\b_\TT$.\qed 
\end{corollary} 
 
Note here that the stalk $\F_\sigma$ of a presheaf $\F$ on $\Delta$ 
coincides with $\F(\sigma)$, since for each point~$\sigma$ of the fan
space, the basic open set~$\langle\sigma\rangle$ is its smallest open 
neighbourhood. 
 
In the simplicial case, that homomorphism turns out to be an 
isomorphism: 
 
\begin{theorem}\label{isom1.3} 
For a simplicial fan $\Delta$, the homomorphism 
$$ 
  \H\b_\TT \longto \A\b 
$$ 
of graded presheaves is an isomorphism, so $\H\b_\TT$ is a sheaf (and 
actually flabby). 
\end{theorem} 
 
\begin{proof} 
We proceed by induction on the number of cones in the fan~$\Delta$. 
For $\Delta = \{o\}$, the assertion is obvious. For the induction 
step, we choose a maximal cone $\sigma \in \Delta$ and consider the 
Mayer-Vietoris sequences associated to $\Lambda := \Delta \setminus 
\{ \sigma \}$ and $\langle\sigma\rangle$, both for $\H\b_\TT$ and 
for $\A\b$. It suffices to prove that in the commutative~diagram 
\[ 
\begin{array}{ccccccccc} 
0 &\!\! \longto \!\! & \H^{2q}_\TT(\Delta) & \!\! \longto\!\! 
& \H^{2q}_\TT(\Lambda) \oplus \H^{2q}_\TT(\sigma) & \!\! \buildrel 
\alpha \over \longto\!\! & \H^{2q}_\TT(\Lambda \cap \sigma) 
& \!\!\longto\!\! & \!\! \H^{2q+1}_\TT(\Delta)\\[3pt] 
& & \mimapdown{\cong} & & 
\mimapdown{\cong} & & \mimapdown{\cong} & & \\[3pt] 
0 & \!\! 
\longto\!\! & \A^{2q}(\Delta) & \!\! \longto\!\! & \A^{2q}(\Lambda) 
\oplus \A^{2q}(\sigma) & \!\! \buildrel \beta \over \longto\!\! & 
\A^{2q}(\Lambda \cap \sigma) & \!\!\longto\!\! & 0 
\end{array} 
\] 
obtained from Corollary~1.2, the rows are exact, the vertical arrows 
are isomorphisms, and $\H^{2q+1}_{\TT}(\Delta)$ vanishes. Applying 
the induction hypothesis to the fans $\Lambda$ and $\Lambda \cap 
\langle\sigma\rangle$, and Lemma~1.1 to~$\langle\sigma\rangle$, we see 
that the assertion holds for the second and the third arrow, 
and we obtain the leading~$0$ in the 
upper row. Furthermore, since the fan $\Delta$ is simplicial, it is 
not difficult to see that the sheaf~$\A\b$ is flabby; 
hence, the map~$\beta$, and thus also~$\alpha$, is surjective. This 
proves $\H^{2q+1}_{\TT}(\Delta) = 0$, as, by induction hypothesis, we 
know that $\H^{2q+1}_{\TT}(\Lambda) \oplus \H^{2q+1}_{\TT}(\sigma)$ 
vanishes. By the Five Lemma, the first vertical arrow is an 
isomorphism as well, thus proving our~claim.
\end{proof} 
 
In the non-simplicial case, we still have some partial results which 
are proved by similar arguments: 
 
\begin{remark} \label{rem1.4} (i) If a fan $\Delta$ can be successively 
built up from cones $\sigma_1, \dots, \sigma_r$ such that $\sigma_{i+1}$ 
intersects $\sigma_1 \cup \dots \cup \sigma_i$ in a single proper face, 
then the homomorphism $\H\b_\TT \to \A\b$ is an isomorphism. 
 
(ii) For an arbitrary toric variety, we have isomorphisms of sheaves 
$$
  \H^q_\TT \cong \A^q \quad\hbox{in all degrees $q \le 2$;} 
$$ 
in particular, $\H^2_\TT$ is a sheaf, and there is an isomorphism 
$$ 
  H^2_\TT(X_{\Delta}) \cong \A^2(\Delta)\; . 
$$ 
\end{remark} 
 
On the other hand, this does not carry over to degree~$q=3$, as 
follows from the next example: 
 
\begin{example} \label{exam1.5} There is a three-dimensional toric 
variety~$X$ with 
$ 
  H^3_\TT(X) \cong \QQ\, . 
$ 
\end{example} 
 
Since we know that $\H^3_\TT(\sigma) \cong \A^3(\sigma)$ vanishes on 
each basic open set, the associated sheaf is the zero sheaf. Hence the 
example implies the following statement: 
 
\begin{corollary} \label{cor1.6} In the non-simplicial case, the presheaf 
$\H\b_\TT$ is in general not a sheaf. 
\end{corollary} 
 
\begin{prexao} 
We consider the fan~$\Delta$ generated by the four ``vertical'' facets 
of a cube centred at the origin, and write it in the form $\Delta_{1} 
\cup \Delta_{2}$, where~$\Delta_{1}$ is generated by two adjacent 
facets, and~$\Delta_{2}$ is generated by the other two. According to 
Remark~1.4, i), we have isomorphisms $\H\b_\TT(\Delta_{j}) \cong 
\A\b (\Delta_{j})$ and $\H\b_\TT(\Delta_{1} \cap \Delta_{2}) \cong 
\A\b (\Delta_{1} \cap \Delta_{2})$. We now show that the homomorphism 
$$ 
  \A^2(\Delta_{1}) \oplus \A^2(\Delta_{2}) \longrightarrow 
  \A^2(\Delta_{1} \cap \Delta_{2}) 
$$ 
in the appropriate Mayer-Vietoris sequence is not surjective. In fact, 
its cokernel is one-dimen\-sio\-nal: The fan $\Delta_{1} \cap 
\Delta_{2}$ is the union $\langle \tau_{1} \rangle \cup \langle 
\tau_{2} \rangle$ of the subfans generated by the two opposite 
two-dimensional ``vertical'' cones $\tau_{j}$ that are spanned by the 
``outer'' vertical edges of the two adjacent facets. The vector space 
$\A^2(\Delta_{1} \cap \Delta_{2}) = A^2_{\tau_1} \oplus A^2_{\tau_2}$ 
is four-dimensional. The restriction homomorphisms from $\A^2(\Delta_{j})$ 
to $\A^2(\Delta_{1} \cap \Delta_{2})$ both map onto the threee-dimensional 
subspace consisting of all piecewise linear functions for which the 
differences between the values at the top vertex and at the bottom 
vertex on both edges agree. Applying that $\H\b_\TT(\Lambda) = 
\A\b(\Lambda)$ holds for $\Lambda = \Delta_{j}$ and for $\Lambda = 
\Delta_{1} \cap \Delta_{2}$, we get $H^3_\TT(X_{\Delta}) \cong \QQ$. 
\hfill $\Box$ 
\end{prexao} 
 
\typeout{1.B. Toric Line Bundles etc.} 
\section*{1.B. Toric Line Bundles and their Equivariant Chern Class} 
 
\noindent A line bundle $L \to X$ on a toric variety~$X$ is called a 
{\it toric line bundle\/} if there is a $\TT$-action by bundle 
automorphisms on the total space such that the bundle projection is 
equivariant. We obtain an induced line bundle $L_{\TT} := E\TT 
\times_{\TT} L$ on $X_\TT$, whose Chern class $c_1(L_{\TT}) \in 
H^2(X_\TT) = H^2_\TT(X)$ is called the {\it equivariant Chern class\/} 
of~$L$ on~$X$ and is denoted with $c_1^\TT(L)$. This class can be 
considered as a lift of the ``usual'' Chern class $c_1(L) \in H^2(X)$ 
to $H^2_\TT(X) = H^2(X_\TT)$: Each fibre map $X \hookrightarrow 
X_\TT$ of the bundle $X_\TT \to B\TT$ induces the same ``edge 
homomorphism'' $H^2_\TT(X) \to H^2(X)$, and that edge homomorphism 
maps $c_1^\TT(L)$ to $c_1(L)$. 
 
Now let $\Delta$ be a {\em purely $n$-dimensional fan} (i.e., a fan 
generated by its $n$-dimensional cones). For use in section~6, we have 
to determine the $\Delta$-piecewise linear function $\psi := \psi_{L} 
\in \A^2(\Delta)$ corresponding to the equivariant Chern class 
$c_1^\TT(L) \in H^2_\TT(X)$ of a toric line bundle~$L$ on~$X$ with 
respect to the isomorphism $H^2_\TT(X_\Delta) \to \A^2(\Delta)$. It is 
convenient to introduce the notation $\check\TT := \TT \times \CC^*$, 
$\check N := N \oplus \ZZ$, and $\check V := V \oplus \QQ$. 
 
The total space of~$L$ is in a natural way a toric variety with acting 
torus $\check\TT$, where the second factor acts on the fibres by scalar 
multiplication. The associated principal $\CC^*$-bundle $L_0$, obtained 
from $L$ by removing the zero section, is an invariant open subset. Let 
us first describe the fan $\check\Delta_0$ in $\check V$ corresponding 
to~$L_{0}$. The fan is determined by the following properties: The 
projection $p \colon \check V \to V$ maps $|\check\Delta_0|$ 
homeomorphically onto $|\Delta|$, inducing a bijection between the 
cones in $\check\Delta_0$ and $\Delta$. In order to assure local 
triviality of the projection, we have to require the equality 
$$ 
  p(\check N_{\check\sigma}) = N_{p(\check\sigma)} 
$$ 
for each cone $\check\sigma$ in $\check\Delta_0$. Then $L$ is the 
$\check\TT$-toric variety given by the fan 
$$ 
  \check\Delta := \check\Delta_0 \cup \{ \check\sigma + \rho \;;\; 
  \check\sigma \in \check\Delta_0 \} \;,\quad \rho := 
  \QQ_{\ge 0} \mal (0_{V},1_{\QQ})
$$ 
spanned by $\check\Delta_0$ and the ``vertical'' ray $\rho$ in $\check V$. 
The support $|\check\Delta_0|$ is the graph of a piecewise linear function 
$\psi := \psi_{L} \colon |\Delta| \to \QQ$ taking 
integral values at lattice points. Vice versa, such a function clearly 
determines a fan~$\check\Delta_0$ of the above type. We note 
that~$\psi_{L}$ is the composition of the $\Delta$-piecewise linear 
homeomorphism $(p|_{\check\Delta_0})^{-1} \colon |\Delta| \to
|\check\Delta_0|$ and the linear projection $\check V \to \QQ$ onto the 
second factor. 
 
\begin{remark} \label{rem1.7} 
If, as above, the projection $p \colon \check V \to V$ maps 
$|\check\Delta_0|$ homeomorphically onto $|\Delta|$, but if the condition 
$p(\check N_{\check\sigma}) = N_{p(\check\sigma)}$ is not satisfied, then 
$L_0$ is called a {\it toric Seifert bundle\/}. Replacing the lattice 
$\check N = N \times \ZZ$ by $N \times \frac{1}{m}\ZZ$ for a suitable 
$m\in \NN_{>0}$, we see that the above condition holds for that new 
lattice. On the level of tori and toric varieties, that means passing 
from $\TT \times \CC^*$ to $\TT \times \CC^*/C_m \cong \TT \times \CC^*$ 
and from $L_0$ to $L_0/C_m$, respectively, where the group $C_m \subset 
\CC^*$ of $m$-th roots of unity acts on $L_0$ as subgroup of the second 
factor in $\TT \times \CC^*$. 
\end{remark} 
 
\begin{lemma} \label{lem1.8} 
For a purely $n$-dimensional fan~$\Delta$, the function $\psi_{L} \in 
A^2(\Delta)$ is the image of the equivariant Chern class $c_1^\TT(L)$ 
of $L$ under the isomorphism $H^2_\TT(X_\Delta) \to \A^2 (\Delta)$ 
from {\/\rm Remark 1.4, (ii)}. 
\end{lemma} 
 
\begin{proof} 
For a cone $\sigma \in \Delta^n$, let $\psi_\sigma \in M$ be the 
character which coincides with $\psi_L$ on $\sigma \cap N$. Since the 
map $\H^2_\TT \to \A^2$ is an isomorphism of sheaves, it suffices to 
show that the ``local'' equivariant Chern class $c_1^\TT(L|_{X_\sigma}) 
\in H^2_\TT(X_\sigma)$ is mapped onto $\psi_\sigma \in M \cong 
\A^2(\sigma)$. Observing that the inclusion $x_\sigma \hookrightarrow 
X_\sigma$ of the fixed point $x_\sigma \in X_\sigma$ induces an 
isomorphism $H^2_\TT(X_\sigma) \buildrel \cong \over \longto 
H^2_\TT(x_\sigma)$, we may further restrict our attention 
to the fibre $L_{x_\sigma}$ of $L$ over~$x_\sigma$. As a $\TT$-module, 
this fibre is nothing but $\CC_{\psi_\sigma}$, and the character 
$\psi_\sigma \in M \cong H^2(B\TT) \cong H^2_{\TT}(x_{\sigma})$ is the 
Chern class of that bundle. This completes the proof. 
\end{proof}

\typeout{2. Equivariant Intersection Cohomology Sheaf} 
\section{The Equivariant Intersection Cohomology Sheaf} 
\setcounter{theorem}{0} 
 
\noindent For a non-simplicial fan~$\Delta$, the equivariant cohomology 
presheaf $\H\b_\TT$ is no longer a sheaf (see Example~1.5), so 
Theorem~1.3 fails to be true in the general case. The situation is much 
better behaved for {\it intersection\/} cohomology, though we do not 
have such a nice combinatorial interpretation as is given by the 
Stanley-Reisner ring in the simplicial case. \par 
 
First let us recall how to describe equivariant intersection 
cohomology: Following the approach of F.~Kirwan\footnote{See \cite{Ki}, 
formula (2.12) and the surrounding text; for a more ``sophisticated'' 
approach, see \cite{Bry} and \cite{Jo}.} and using the notation 
of section~1, one defines the $q$-th {\it(rational) equivariant 
intersection cohomology group\/} of $X$ as the limit 
$$ 
  IH^q_\TT(X) := IH^q_\TT(X; \QQ) := 
  \lim_{m \to \infty} IH^q(X_{\TT,m}; \QQ)
$$ 
and sets 
$$ 
  IH\b_\TT(X):= \bigoplus_{q=0}^\infty IH^q_\TT(X)\, . 
$$ 
This construction provides a presheaf 
$$ 
  \IH\b_\TT 
  \colon \Lambda \mapsto \IH\b_\TT (\Lambda) 
  := IH\b_\TT(X_\Lambda) \quad\hbox{for $\Lambda \preceq \Delta$} 
$$ 
on the fan space~$\Delta$. In order to prove that it is in fact a 
sheaf, we verify the following basic result: 
 
\begin{vanlemma} \label{vanlem} 
The equivariant intersection cohomology $IH_\TT^q(X)$ of a toric variety 
$X$ vanishes in odd degrees~$q$. 
\end{vanlemma} 
 
\begin{proof} Let $\hat\Delta$ be a simplicial refinement of the 
defining fan~$\Delta$ for~$X$, and denote with $\hat X \to X$ the 
corresponding equivariant $\QQ$-resolution of singularities. By 
theorem~1.3, the assertion holds for~$\hat X$, since then 
$IH\b_\TT(\hat X) \cong H\b_\TT(\hat X) \cong \A\b_{\hat\Delta}$. 
By the equivariant version of the Decomposition Theorem of Beilinson, 
Bernstein, 
Deligne, and Gabber as stated in \cite[p.~394]{Ki} (see also
\cite[5.3]{BeLu} 
or \cite{BreLu}), we may interpret $IH^q_{\TT}(X)$ as a subspace of 
$IH^q_{\TT}(\hat 
X)$. That proves the assertion. \end{proof} 
 
\begin{theorem} \label{sheaf2.2} 
The presheaf $\IH\b_\TT$ 
on $\Delta$ is a sheaf of $\A\b$-modules. 
\end{theorem} 
 
\begin{proof} Since there are only finitely many open subsets
  in~$\Delta$, 
it suffices to verify the sheaf axioms for two open
subsets~$\Lambda_1, 
\Lambda_2 \preceq \Delta$. We thus have to prove the exactness of the
sequence 
$$ 
0 \longto \IH_\TT^q(\Lambda_1 \cup \Lambda_2) \longto
\IH_\TT^q(\Lambda_1) 
\oplus \IH_\TT^q(\Lambda_2) \longto \IH_\TT^q(\Lambda_1 \cap
\Lambda_2) 
\;. 
$$ 
That follows from the Lemma~2.1: The exactness is obvious 
if $q$ is odd; for even~$q$, the vector space $\IH_\TT^{q-1}(\Lambda_1 
\cap \Lambda_2)$ vanishes and thus, the sequence is part of 
the exact Mayer-Vietoris sequence for~$IH_\TT\b$. 
 
As a consequence, $\IH\b_\TT$ is a sheaf of $\A\b$-modules, since
$\A\b$ 
is the associated sheaf to the presheaf $\H\b_\TT$, and each 
$IH\b_\TT(X_\Lambda)$ 
is an $H\b_\TT(X_\Lambda)$-module. \end{proof} 
 
\typeout{3. Minimal Extension Sheaves} 
\section{Minimal Extension Sheaves} 
\setcounter{theorem}{0} 
 
\noindent We now proceed towards (re-)stating and verifying the three 
properties that actually characterize the sheaf $\IH\b_\TT$ up to 
isomorphism.

For the ease of notation, we write $E\b_\Lambda$ for the module 
$\E\b(\Lambda)$ of sections of a sheaf $\E\b$ over a subfan $\Lambda 
\preceq \Delta$; in particular, we do so for the residue class module 
$\overline{E\b_\Lambda} := E\b_\Lambda/(\mm \mal E\b_\Lambda) \cong 
\QQ\b \otimes_{A\b} E\b_\Lambda$, see~(0.2). Using this notation, we 
restate the properties (N) ({\it Normalization\/}), (PF) ({\it 
Pointwise Freeness\/}), and (LME) ({\it Local Minimal Extension 
$\bmod\; \mm$\/}) from the introduction.
 
\begin{definition} \label{mes} 
A sheaf $\E\b$ of graded $\A\b$-modules on the fan $\Delta$ is called 
a {\it minimal extension sheaf\/}(of~$\QQ\b$) if it satisfies the 
following conditions: 
 
\begin{description} 
\item[(N)] One has $E\b_o \cong A\b_o = \QQ\b$ for the zero cone~$o$. 
\item[(PF)] For each cone $\sigma \in \Delta$, the module~$E\b_{\sigma}$ 
  is {\it free} over~$A\b_{\sigma}$. 
\item[(LME)] For each cone $\sigma \in \Delta$, the 
  restriction mapping $\varphi_\sigma \colon E\b_{\sigma} \to 
  E\b _{\partial \sigma}$ induces an isomorphism 
$$ 
  \overline \varphi_\sigma \colon \overline{E}\b_{\sigma} \buildrel 
  \cong \over \longrightarrow \overline{E}\b_{\partial\sigma} 
$$ 
  of graded vector spaces. 
\end{description} 
\end{definition} 
 
Condition~(LME) implies that $\E\b$ is minimal in the set of all flabby 
sheaves of graded $\A\b$-modules that satisfy the conditions~(N) 
and~(PF), cf.\ Remark~3.3, whence the name ``minimal extension sheaf''. 
-- Furthermore, let us note that on a simplicial subfan~$\Lambda$, the 
restriction~$\E\b|_{\Lambda}$ of such a sheaf is isomorphic 
to~$\A\b|_{\Lambda}$, so it is locally free of rank~one. In the 
non-simplicial case, however, the rank of the stalks of $\E\b$ is not 
constant, so $\E\b$ can not be a locally free $\A\b$-module. 
 
We now show how to construct such a sheaf on an arbitrary fan space. 
 
\begin{proposition} \label{prop3.1} 
On every fan $\Delta$, there exists a minimal extension sheaf $\E\b$, 
and it is unique up to isomorphism. 
\end{proposition} 
 
\begin{proof} We define the sheaf $\E\b$ inductively on the $k$-skeleton 
subfans 
$$ 
  \Delta^{\le k} := \bigcup_{j \le k} \Delta^j \;,\quad 
  \Delta^j := \{\sigma \in \Delta \;;\; \dim\sigma = j\}\;, 
$$ 
starting with $E_o \b := \QQ\b$ for $k=0$. Suppose that for some 
$k > 0$, the sheaf~$\E\b$ has already been constructed on 
$\Delta^{\le k-1}$, so in particular, for each $k$-dimensional cone 
$\sigma$, the module $E\b_{\partial\sigma}$ is given. It thus suffices 
to define $E\b_{\sigma}$ together with a restriction homomorphism 
$E\b_{\sigma} \to E\b_{\partial \sigma}$ inducing an isomorphism on the 
quotients modulo~$\mm$. This is achieved by setting 
$$ 
E\b_{\sigma} := A\b_{\sigma} \otimes_{\QQ} \overline E\b_{\partial
  \sigma} 
$$ 
with the restriction map being induced by some $\QQ\b$-linear section 
$s: \overline E\b_{\partial\sigma} \to E\b_{\partial \sigma}$ of the
residue 
class map $E\b_{\partial \sigma} \to \overline E\b_{\partial
  \sigma}$. 
 
The unicity is proved similarly in an inductive manner; we refer to 
our companion article \cite{BBFK$_{2}$} for details. \end{proof} 
 
\begin{remark} \label{rem3.2} 
A minimal extension sheaf $\E\b$ is flabby and vanishes in odd
degrees. 
\end{remark} 
 
\begin{proof} Since the fan space~$\Delta$ is covered by finitely many 
affine fans, it suffices to prove that for each cone $\sigma \in 
\Delta$, the restriction map $\varphi_\sigma$ is surjective. Using the 
results on graded modules recalled in the ``Preliminaries'', that is a 
consequence of condition~(LME). -- The vanishing of~$\E^{2q+1}$ follows 
immediately from the same condition, since $A\b_\sigma$ and thus $E\b_o 
\cong A\b_o = \QQ\b$ ``live'' only in even degrees. 
\end{proof} 
 
The properties of minimal extension sheaves are investigated in 
\cite{BBFK$_{2}$}; let us quote just one result: 
\begin{quote} 
  {\it The sheaf $\A\b$ is a minimal extension sheaf if and only 
  if the fan~$\Delta$ is simplicial.} 
\end{quote} 
Our aim here is to show that $\IH\b_\TT$ represents this unique 
isomorphism class of minimal extension sheaves.
 
\begin{mntheorem} \label{mntheo} 
The equivariant intersection cohomology sheaf $\IH\b_\TT$ on $\Delta$ 
is a minimal extension sheaf. 
\end{mntheorem} 
 
\begin{proof} The Normalization property is obviously satisfied, since 
we have 
$$
  \IH\b_\TT(o) = IH\b_\TT(\TT) \cong \QQ\b\, . 
$$ 
The Pointwise Freeness condition will be verified in Corollary~4.5, and 
the Local Minimal Extension requirement, in Proposition~5.1. 
\end{proof} 
 
\typeout{4. Equivariantly Formal Toric Varieties etc.} 
\section{Equivariantly Formal Toric Varieties and Pointwise Freeness 
{of~$\IH\b_\TT$}.} 
\setcounter{theorem}{0} 
 
\noindent We now proceed to verify the condition 
\begin{enumerate} 
  \item[(PF)] {\it For each cone~$\sigma$, the $A\b_{\sigma}$-module 
  $IH\b_{\TT}(X_{\sigma})$ is free\/.} 
\end{enumerate} 
This follows immediately from the isomorphism 
$$ 
  IH\b_{\TT}(X_{\sigma}) \cong A\b_{\sigma} \otimes IH\b(Z_{\sigma}) 
$$ 
obtained in the proof of Corollary~4.5 below, where $Z_{\sigma}$ is the 
contractible affine $\TT_{\sigma}$-toric variety occuring in the orbit 
splitting $X_{\sigma} \cong {\TT\,}' \times Z_{\sigma}$ as in~(0.3). 
The crucial point is that~$Z_{\sigma}$ satisfies the conditions -- of 
course, with respect to the acting torus~$\TT_{\sigma}$ -- stated in 
the following result, a more general version of which can be found in 
\cite{GoKoMPh}. 
 
\begin{lemma} \label{lem4.1} 
For a toric variety $X$, the following statements are equivalent: 
\begin{enumerate} 
\item[i)] The K\"unneth formula $IH\b_\TT(X) \cong H\b(B\TT) \otimes 
  IH\b(X)$ holds (which, as above, implies that $IH_\TT\b(X)$ is a 
  free $A\b$-module). 
\item[ii)] With $\overline{IH\b_\TT(X)} := 
  IH\b_\TT(X)/\big(\mm \mal IH\b_\TT(X)\big)$, each inclusion $X 
  \hookrightarrow X_\TT$ as a fibre induces an isomorphism 
$$ 
  \overline{IH\b_\TT(X)} \;\;\cong\;\; IH\b(X) \leqno{(4.1)} 
$$ 
  of graded vector spaces. 
\item[iii)] The non-equivariant intersection cohomology $IH^q(X)$ 
  vanishes in odd degrees~$q$. 
\end{enumerate} 
\end{lemma} 
 
\begin{proof} Since condition~i) says that for intersection 
cohomology, $X_\TT$ behaves like the product $X \times B\TT$, the 
implication ``i) $\Rightarrow$ ii)'' is obvious, and ``ii) 
$\Rightarrow$ iii)'' follows immediately from the Vanishing 
Lemma~2.1. For the implication ``iii) 
$\Rightarrow$ i)'', we observe that the assumption implies the 
degeneration of the intersection cohomology Leray-Serre spectral 
sequence 
$$ 
  E_2^{p,q} = H^p(B\TT) \otimes IH^q(X) \Rightarrow IH^{p+q}_\TT(X) 
$$ 
associated to the fibering $X_\TT \to B\TT$ at the $E_{2}$-level: 
Since $H^p(B\TT)$ vanishes for odd $p$, the spectral terms $E_2^{p,q}$ 
vanish for odd total degrees $k = p+q$, and consequently, the 
differentials $d_2^{p,q} \colon E_2^{p,q} \to E_2^{p+2,q-1}$ are trivial. 
By induction on~$r$, that holds also for every $r \ge 2$. 
\end{proof} 
 
Thus, for toric varieties satisfying these properties, the equivariant 
and the non-equi\-variant theory determine each other in a simple way. 
Following \cite{GoKoMPh}, we use the following terminology: 
 
\begin{definition} \label{equifor} 
A toric variety $X$ -- and its defining fan -- are called 
{\em equivariantly formal\/} (for intersection cohomology), or 
$IH_\TT$-formal for short, if and only if $X$ satisfies one -- 
and hence all -- of the above three conditions. 
\end{definition} 
 
We point out that there is an analogous notion for ``ordinary'' 
cohomology, and that a toric variety may be $IH_\TT$-formal, but not 
$H_\TT$-formal (e.g., a compact toric threefold with $b_3 \ne 0$). 
As we are mainly dealing with intersection cohomology, however, we use 
``equivariantly formal'' only in the sense of ``$IH_\TT$-formal''. 
 
The most important cases when toric varieties are equivariantly formal 
are the compact case and the contractible affine case. From the theorem 
of Jurkiewicz and Danilov cited in the introduction, it follows that a 
rationally smooth compact toric variety is equivariantly formal. The 
proof in the compact case is now an easy consequence of the famous 
Decomposition Theorem: 
 
\begin{proposition} \label{compa} 
A compact toric variety $X$ is equivariantly formal. 
\end{proposition} 
 
\begin{proof} Let $\pi: \hat X \to X$ denote a toric $\QQ$-resolution 
as in the proof of the Vanishing Lemma~2.1. By the ``classical'' (i.e, 
non-equivariant) Decomposition Theorem of Beilinson, Bernstein, 
Deligne and Gabber, we know that $IH\b(X)$ is isomorphic to a direct 
summand of $H\b(\hat X)$, and, according to the theorem of 
Jur\-kie\-wicz and Danilov, that module vanishes in odd degrees. 
\end{proof} 
 
We now state the analogous result for the contractible affine case 
where it is considerably more difficult to handle: 
 
\begin{theorem} \label{affine} 
A contractible affine toric variety $X$ is equivariantly formal. 
\end{theorem} 

This result immediately yields the ``Pointwise Freeness'' 
property~(PF).

\begin{corollary} \label{cor4.4} 
For each cone $\sigma \in \Delta$, the $A\b_{\sigma}$-module 
$\IH\b_\TT(X_\sigma)$ is free. 
\end{corollary} 
 
\begin{procoo} 
The affine orbit splitting $X_{\sigma} \cong {\TT\,}' \times Z_\sigma$ 
as in~(0.3) yields the isomorphism 
$$ 
IH\b_\TT(X_\sigma) \cong IH\b_{{\TT\,}' \times \TT_{\sigma}}({\TT\,}' 
\times Z_\sigma) \cong IH\b_{\TT_\sigma}(Z_\sigma) \;.\leqno{(4.2)} 
$$ 
By the theorem, the contractible affine $\TT_{\sigma}$-toric 
variety~$Z_{\sigma}$ is $IH\b_{\TT_{\sigma}}$-formal, i.e., we have 
isomorphisms 
$$ 
  IH\b_{\TT_\sigma}(Z_\sigma) \;\cong \; 
  H\b (B\TT_{\sigma}) \otimes IH\b(Z_\sigma) 
  \;=\; A\b_\sigma \otimes IH\b (Z_\sigma)\; , 
$$ 
showing that $IH\b_{\TT_\sigma}(Z_\sigma)$ is a free 
$A\b_\sigma$-module and thus proving the corollary. 
\hfill$\Box$ 
\end{procoo} 
 
\begin{protheo} Writing $X := X_{\sigma}$ for ease of notation, we 
first notice that $X$ is a distinguished neighbourhood of its (unique) 
fixed point~$x$. According to the attachment condition (see 
\cite[V.4.2]{Bo}), we thus have an isomorphism 
$$ 
  IH\b(X) \cong \tau_{<n} IH\b(X^*)\, ,\leqno{(4.3)} 
$$ 
where $X^* := X \setminus \{x\}$ is the ``punctured'' toric variety 
obtained by removing the fixed point. Hence, it suffices to prove that 
$IH^q(X^{*})$ vanishes in odd degrees $q < n$. The basic idea to reach 
that aim is to pass from~$X^*$ to a projective toric variety~$Y$ 
having~$X$ as its ``affine cone'' and then to compare $IH\b(X^*)$ 
with~$IH\b(Y)$, keeping in mind that $Y$ is equivariantly formal. 
 
Such a projective toric variety~$Y$ is obtained as quotient of~$X^{*}$ 
modulo the action of any one parameter subgroup of~$\TT$ having the 
following property: The orbits of the induced $\CC^*$-action on~$X$ 
have the fixed point~$x$ as common ``source'', so an equivariant 
contraction is provided by $(t,x) \mapsto t \mal x$ for $t \to 0$, 
where~$t \in\; ]0,1] \subset \CC^{*}$ is the parameter. A one parameter 
subgroup satisfies these conditions if and only if the representing 
lattice vector $\alpha \in N = \Hom(\CC^*, \TT)$ lies in the relative 
interior of~$\sigma$. If, in addition, the lattice vector~$\alpha$ is 
primitive, then the induced $\CC^*$-action on~$X^{*}$ is free, possibly 
up to some finite isotropy along lower-dimensional $\TT$-orbits. 
Let~$F \subset \alpha(\CC^*)$ be the finite subgroup generated by these 
isotropy groups. Then the equivariant mapping $X^{*} \to Y$ factors 
through the quotient $X^{*}/F$, making $X^{*}/F \to Y$ a principal 
$\CC^{*}$-bundle (cf.\ Remark~1.7). Replacing~$X^{*}$  with $X^{*}/F$ 
does not change the intersection cohomology: In fact, there  is an 
identification 
$$ 
  IH\b(X^*/F) \cong IH\b(X^*)^F = IH\b(X^*) \leqno{(4.4)} 
$$ 
(see \cite[Lemma 2.12]{Ki$_{2}$}), the equality holding since the action 
of $F$ on $IH\b(X^*)$ is trivial: It is induced from a -- necessarily 
trivial -- action of the connected group $\TT$ on the rational vector 
space $IH^*(X^*)$. 
 
The quotient $Y$ is a toric variety for the torus $\overline \TT := 
\TT / \alpha(\CC^*)$, whose lattice of one parameter subgroups is 
$\overline{N} := N/(\ZZ \mal \alpha)$. The defining fan for~$Y$ is 
the image 
$$ 
 \Phi := p(\partial\sigma) := \{p(\tau) \;;\; \tau \prec \sigma\} 
$$ 
of the boundary fan $\partial\sigma$ under the quotient mapping 
$p \colon V \to W$ from $V = N_{\QQ}$ onto the quotient vector space 
$W := \overline{N}_{\QQ} := N_{\QQ}/(\QQ \mal \alpha)$. 
 
We now proceed to proving that $\tau_{<n}IH\b(X^*)$ vanishes in odd 
degrees, thus showing that $X$ is equivariantly formal. By the 
identification (4.4), we may assume without loss of generality 
that~$F$ is the trivial subgroup and thus, that $X^* \to Y$ is a 
principal $\CC^*$-bundle over the projective toric variety $Y$. We 
next consider the toric line bundle $L \to Y$ obtained from $X^* 
\to Y$ by adding a zero section opposite to the fixed point of $X$ 
(``section at infinity'') to that $\CC^*$-bundle: Its total 
space~$L$ is the toric variety associated to the fan 
$$ 
  \Sigma := \partial\sigma \cup 
  \{ \tau + \QQ_{\ge 0}\mal(-\alpha)\;;\; \tau \prec \sigma \} 
$$ 
in $V = N_{\QQ}$, with the projection $N \to \overline N$ inducing a 
mapping of fans $\Sigma \to \Phi$ and thus, a toric morphism 
$L \to Y$. The line bundle~$L \to Y$ is ample (see, e.g., 
\cite[2.12, p.~82]{Oda$_1$}), so a suitable tensor power is very 
ample. If~$L$ is already very ample itself, then the one point 
compactification of~$L$ is the projective cone over~$Y$ with respect 
to the projective embedding determined by the sections of~$L$, the 
corresponding complete fan being $\Sigma \cup \{ \sigma \}$. 
 
Let us now look at the following commutative diagram whose top row is 
part of the long exact $IH\b$-sequence associated to the pair
$(L,X^*)$: 
\[ 
\begin{array}{ccccccccc} 
\dots \longto\!\! & IH^{q-1}(X^*) & 
\!\!\longto\!\! & IH^q(L,X^*) & \!\!\longto\!\! & IH^q(L) & 
\!\!\longto\!\! & IH^q(X^*) & \!\!\longto \dots\;{} \\[3pt] 
& & & \mimapup{\cong} & & \mimapup{\cong} & & & \\[3pt]
& & & IH^{q-2}(Y) & \!\!\buildrel \lambda\over\longto\!\! & 
IH^q(Y) & & & \; 
\end{array} 
\] 
Here the first vertical isomorphism is the Thom isomorphism for the 
line bundle, the second one is induced by the bundle projection, and 
$$
  \lambda \colon IH\b(Y)[-2] \to  IH\b(Y)
$$ 
is the homomorphism given by cup product multiplication with 
$c_{1}(L) \in H^2(Y)$. The resulting long exact 
sequence 
$$ 
 \dots \longto IH^{q-1}(X^*) \longto IH^{q-2}(Y) 
 \buildrel \lambda \over \longto IH^q(Y) \longto IH^q(X^*) \longto 
 \dots 
$$ 
is the Gysin sequence associated to the bundle $X^{*} \to Y$. 
 
Using the fact that~$Y$ is equivariantly formal and thus has vanishing 
intersection cohomology in odd degrees, the Gysin sequence decomposes 
into shorter exact sequences 
$$ 
0 \longto IH^{q-1}(X^*) \longto IH^{q-2}(Y) 
\buildrel \lambda \over \longto IH^q(Y) \longto IH^q(X^*) \longto 0 
$$ 
if $q$ is even. Hence, it suffices to show that $\lambda$ is injective 
for $q \le n$. 
 
Since $L$ is ample, we may apply the hard Lefschetz theorem for 
intersection cohomology: a suitable tensor power $L^{\otimes m}$ of 
$L$ is very ample, and since $c_1(L^{\otimes m})= m c_1(L)$, the 
assertion of the theorem holds for $L$ as well. Hence, we know that 
$$ 
\lambda^k \colon IH^{n-1-k}(Y) \longto IH^{n-1+k}(Y) 
$$ 
is an isomorphism for every $k \ge 0\,$; as a consequence, we see that 
$$ \lambda \colon IH^{q-2}(Y) \longto IH^q(Y) $$ 
is injective for $q \le n$ and surjective for $q \ge n$. 
\hfill $\Box$ 
\end{protheo} 
 
For later use, we note that the properties of~$\lambda$ yield 
an isomorphism 
$$ 
  \coker\,\lambda \cong \tau_{<n}IH\b(X^*) \,.\leqno{(4.5)} 
$$ 
 
\typeout{5. Local Minimal Extension Property etc.} 
\section{Local Minimal Extension Property of $\IH\b_\TT$.} 
\setcounter{theorem}{0} 
 
\noindent To complete the proof of our Main Theorem stating that the 
sheaf $\IH\b_{\TT}$ is a minimal extension sheaf, we have to verify 
the condition~(LME). We restate it as follows: 
 
\begin{proposition} \label{prop5.1} 
For each cone~$\sigma$, the restriction mapping 
$$ 
  IH\b_{\TT}(X_{\sigma}) \to IH\b_{\TT}(X_{\partial\sigma}) 
$$ 
induces an isomorphism between the quotients modulo the maximal 
ideal~$\mm$ of~$A\b$. 
\end{proposition} 
 
\begin{proof} By the ``relative affine orbit splitting'' 
$(X_\sigma, X_{\partial\sigma}) \cong \TT' \times 
(Z_\sigma, Z_{\partial\sigma})$ and the formula (4.2) as well as its 
analogue $IH\b_{\TT}(X_{\partial\sigma}) \cong 
IH\b_{\TT_{\sigma}}(Z_{\partial\sigma})$, we see that it is sufficient 
to consider the case of an $n$-dimensional cone~$\sigma$. We use the 
same notations as in the previous section; in particular, we write $X 
= X_\sigma$ and $X^* = X^*_\sigma$. Furthermore, by the arguments of 
the preceding section, we may replace $\TT$, $X_\sigma$, and 
$X_{\partial \sigma}$ by $\TT/F$, $X_\sigma/F$, and 
$X_{\partial \sigma}/F$, respectively, where $F \subset \TT$ is a 
suitable finite subgroup, without changing the base ring $H\b(B\,\TT) 
\cong H\b\bigl(B(\TT/F)\bigr)$ and the above homomorphism. Hence we 
may assume that $X^* \to Y$ is a principal $\CC^*$-bundle. 
 
First we collect in a big commutative diagram all the objects we have 
to consider: 
\[ 
\begin{array}{ccccc} 
  IH\b_\TT(X) & \longto & IH\b_\TT(X^*) & 
    {\mathrel{\mathop{\longleftarrow }\limits^{\cong}_{\pi_{\TT}}}}
    & IH\b_{\overline \TT}(Y) \\[3pt] 
  \mimapdown{} & & \mimapdown{} & & \mimapdown{} \\[3pt] 
  \overline{IH\b_\TT(X)} & \longto & \overline{IH\b_\TT(X^*)} & 
    \buildrel{\overline\pi_{\TT}}\over \longleftarrow & 
    \overline {IH\b_{\overline \TT}(Y)} \\[3pt] 
  \mimapdown{\cong} & & \bigcap & & \mimapdown{\cong} \\[3pt] 
  IH\b(X) & \longto & IH\b(X^*) & 
    \buildrel{\pi}\over \longleftarrow & IH\b(Y) \\[3pt] 
  & \llap{$\scriptstyle\cong$} \searrow & \bigcup & & \mimapdown{} \\[3pt] 
  & & \tau_{<n}IH\b(X^*) & \buildrel \cong \over \longleftarrow &  
    \coker\, \lambda {\,.}
\end{array} 
\] 
Here, $\lambda$ again denotes the homomorphism given by the cup product 
with the Chern class $c_1(L)$ of the line bundle $L \to Y$ as in the 
previous section. 
 
Once having established the diagram, the proof is achieved as soon as 
we have identified the quotient 
$$ 
  \overline{IH\b_\TT(X^*)} \;:= \; (A\b/\mm) \otimes_{A\b} 
  IH\b_\TT(X^*) \;\cong\; 
  {IH\b_\TT(X^*)} \bigm/ {\mm \mal IH\b_\TT(X^*)} 
$$ 
-- this is just the image of the ``edge homomorphism'' 
$IH\b_\TT(X^*) \to IH\b(X^*)$ relating the equivariant and the 
non-equivariant theory -- with $\tau_{<n}IH\b(X^*)$. 
 
To that end, we consider the right hand side of the diagram, carefully 
keeping track of the different module structures in the top row: 
Whereas $IH\b_\TT(X)$ and $IH\b_\TT(X^*)$ are both modules over $A\b = 
S\b(M_\QQ) \cong H\b(B\TT)$, we have to look at 
$IH\b_{\overline \TT}(Y)$ as a module over the ring $B\b := 
S\b(\overline M_\QQ) \cong H\b(B\overline\TT)$, where $\overline M := 
\overline N^*$. In particular, the ``overlined'' modules in the 
second row are quotients modulo the maximal ideals of the respective 
base rings: For~$X$ and~$X^{*}$, this is the ideal $\mm = \mm_{A} := 
A^{>0}$, whereas for~$Y$, we have to consider 
$$ 
  \overline{IH\b_{\overline \TT}(Y)} 
  \;:=\; (B\b/\mm_B) \otimes_{B\b} IH\b_{\overline \TT}(Y) \;\cong\; 
  {IH\b_{\overline \TT}(Y)} \bigm/ 
  {\mm_{B} \mal IH\b_{\overline \TT}(Y)} 
$$ 
(with $\mm_{B} := B^{>0}$). As the projective $\overline \TT$-toric 
variety~$Y$ is equivariantly formal, we may identify this graded 
$\QQ\b$-module with $IH\b(Y)$. We may now consider~$B\b$ as a subring 
of $A\b$ since $\overline M = \overline N^*$ is canonically isomorphic 
to a submodule of $N^*=M$. Hence, writing the horizontal arrow 
$\overline\pi_{\TT} \colon \overline {IH\b_{\overline \TT}(Y)} \to 
\overline{IH\b_\TT(X^*)}$ as 
$$ 
  {IH\b_{\overline \TT}(Y)} 
  \bigm/ {\mm_{B} \mal IH\b_{\overline\TT}(Y)} \longto
  {IH\b_\TT(X^*)} 
  \bigm/ {\mm_{A} \mal IH\b_\TT(X^*)} \;,
$$ 
we see that it is an epimorphism, being induced by the horizontal 
isomorphism~$\pi_{\TT}$ in the top row. Thus, using the isomorphism in 
the bottom row, we are done if we can show that the kernel of~$\pi$ 
equals the image of the ``hard Lefschetz homomorphism''~$\lambda$. 
 
Before we do that let us give some further remarks on the big diagram: 
The isomorphism 
$$ 
  IH\b_{\overline \TT}(Y) \;
  {\mathrel{\mathop{\longto}\limits^{\cong}_{\pi_{\TT}}}}\;
  IH\b_\TT(X^*) \leqno{(5.1)} 
$$ 
in the top row is obtained as follows: The bundle projection 
$p \colon X^* \to Y$ induces a compatible family of projections 
$$ 
  X^*_{\TT,m} \to Y_{\TT,m} \cong Y_{\overline\TT,m} \times 
  B_{m}\CC^{*} \to Y_{\overline\TT,m} 
$$ 
between the finite-dimensional approximations of $X_{\TT}$ and 
$Y_{\overline\TT}$. As these projections are {\it placid\/} maps, 
there is a (unique) induced homomorphism 
$$
  \pi_{\TT} \colon IH\b_{\overline\TT}(Y) \to IH\b_{\TT}(X^*) 
$$ 
(for a discussion, see \cite[\S4]{GoMPh} or \cite[3.3]{BBFGK}). We now 
note that the bundle is locally trivial and that there is a finite 
open affine covering of~$Y$ by toric subvarieties~$V_i$ such that the 
restricted bundle $U_i := p^{-1}(V_i) \to V_i$ actually is trivial. 
By~(4.2), we thus have isomorphisms $IH\b_\TT(U_i) \cong 
IH\b_{\overline \TT}(V_i)$; gluing these by a Mayer-Vietoris argument 
yields the result. -- We note that, by construction, the 
isomorphism~(5.1) is a morphism of $B\b$-modules; it thus induces an 
isomorphism of the quotients modulo the homogeneous maximal ideal 
$\mm_{B} = B^{>0}$, while $\overline{ IH_\TT\b(X^*)}$ refers to the 
bigger base ring $A\b$ and thus is a quotient of 
$\overline{IH_{ \overline \TT}\b(Y)}$. 
 
We recall why we obtain the other isomorphisms occuring in the diagram: 
By the attachment condition, the restriction mapping $IH\b(X) \to 
IH\b(X^*)$ factors through the inclusion of $\tau_{<n}IH\b(X^*)$, 
inducing the oblique isomorphism. The two vertical isomorphisms follow 
from the fact that both~$X$ and~$Y$ are equivariantly formal. The 
lower horizontal isomorphism has been obtained in~(4.5) in the 
previous section. 
 
We now continue with the proof of Proposition 5.1. First of all, we 
lift the ``hard Lefschetz map'' $\lambda \colon IH\b(Y)[-2] \to 
IH\b(Y)$ to a map in equivariant intersection cohomology: We make the 
line bundle $L \to Y$ a $\overline \TT$-toric line bundle by choosing 
some sublattice~$N_{0}$ complementary to~$\ZZ \mal \alpha$, thus 
providing a direct sum decomposition $N = \ZZ \alpha \oplus N_0$. From 
section~1.B, we recall that the $\overline\TT$-equivariant 
Chern class $c_1^{\overline \TT}(L) \in H^2_{\overline\TT}(Y)$ of~$L$ 
is a lifting of the ``usual'' Chern class $c_1(L) \in H^2(Y)$. It follows 
that the cup product with $c_1^{\overline\TT}(L)$ yields a homomorphism 
$\lambda_{\overline \TT} \colon IH\b_{\overline \TT}(Y)[-2] \to 
IH\b_{\overline \TT}(Y)$ that is a lifting of the mapping $\lambda 
\colon IH\b(Y)[-2] \to IH\b(Y)$ given by the cup product with $c_1(L)$. 
 
We further recall that $\IH\b_{\overline\TT}$ is a sheaf on the defining 
fan~$\Phi$ for~$Y$, and it is a module over the ``structure sheaf'' 
${}_\Phi\A\b$ corresponding to the fan $\Phi$. Using the canonical 
isomorphism $H^2_{\overline\TT}(Y) \to \A^2(\Phi)$, we may identify the 
equivariant Chern class~$c_1^{\overline \TT}(L) \in 
H^2(Y_{\overline \TT})$ with the $\Phi$-piecewise linear function $\psi 
:= \psi_{L} \in A^2(\Phi)$ (see Remark~1.4, (ii), and Lemma~1.8). \par 
 
Now the quotient projection $N \to \overline N$ induces an isomorphism 
$$
  {}_\Phi\A\b(\Phi) \buildrel \cong \over \longto
  {}_\Delta\A\b(\partial \sigma)\,,
$$ 
and the image of $\psi \in \A^2(\Phi)$ coincides with the restriction 
to $\partial \sigma$ of the ``global'' linear form $f \in A^2$, the 
projection $f \colon N_\QQ = \QQ \mal \alpha \oplus (N_0)_{\QQ} \to 
\QQ$ mapping $\alpha$ to $-1$ and having $(N_0)_{\QQ}$ as kernel, cf.\ 
section~1.B. 
 
On the other hand, we have an isomorphism of polynomial rings $A\b = 
B\b[f]$ and thus, the equality $\mm = \mm_A = (\mm_B, f) := \mm_B + A\b 
\mal f$ for the homogeneous maximal ideals. The proof of the assertion 
that the restriction homomorphism $IH\b_{\TT}(X) \to IH\b_{\TT}(X^*)$ 
induces an isomorphism on the quotients with respect to the submodules 
generated by~$\mm_A$ is now obtained as follows: Under the inverse of the 
isomorphism~(5.1), the submodule $\mm_{A} \mal IH\b_\TT(X^*) = 
(\mm_{B},f) \mal IH\b_\TT(X^*)$ is mapped onto $(\mm_{B},\psi) \mal 
IH\b_{\overline \TT}(Y)$. We thus have an isomorphism of quotients 
$$ 
  \frac{IH\b_\TT(X^*)}{\mm_{A} \mal IH\b_\TT(X^*)} 
  \;\cong\; \frac{IH\b_{\overline \TT}(Y)} {(\mm_{B},\psi) \mal
  IH\b_{\overline \TT}(Y)} \;\;.\leqno{(5.2)} 
$$ 
As explained above, the mapping $IH\b_{\overline \TT}(Y)[-2] \to 
IH\b_{\overline \TT}(Y)$ given as multiplication by $\psi \in 
\A^2(\Phi)$ lifts the ``Hard Lefschetz homomorphism'' $\lambda \colon 
IH\b(Y)[-2] \to IH\b(Y)$ to the equivariant theory. Since~$Y$, as a 
projective $\overline \TT$-toric variety, is equivariantly formal, we 
may eventually rewrite the right hand side of the above isomorphism 
(5.2) as follows: 
$$ 
  \frac{IH\b_{\overline \TT}(Y)} {(\mm_{B},\psi) \mal 
  IH\b_{\overline \TT}(Y)} \;\cong\; 
  \frac{IH\b(Y)}{\lambda\big(IH\b_{\ }(Y)\big)} 
  = \coker\,\lambda \;. 
$$ 
The proof is now achieved using the isomorphisms 
(4.5), (4.3), and the fact that the contractible affine toric
variety~$X$ is $IH\b_{\TT}$-formal. 
\end{proof} 
 
\typeout{6. Some Results on Equivariantly Formal Fans} 
\section{Some Results on Equivariantly Formal Fans} 
\setcounter{theorem}{0} 
 
\noindent 
We recall from Lemma 4.1 that in the case of an equivariantly formal 
toric variety, the equivariant intersection cohomology determines the 
``usual'' intersection cohomology in a straightforward way, namely, 
as the quotient by the $A\b$-submodule generated by the homogeneous 
maximal ideal $\mm$ of the ring~$A\b$. Since 
the equivariant intersection cohomology sheaf $\IH\b_\TT$ is a minimal 
extension sheaf, we thus have an isomorphism $IH\b(X_\Delta) \cong 
\overline{\E\b(\Delta)} := \E\b(\Delta)/\mm\mal\E\b(\Delta)$, where 
as usual~$\E\b$ denotes a minimal extension sheaf on~$\Delta$. 

It is convenient to call a (rational) fan~$\Delta$ {\it equivariantly 
formal} if the toric variety~$X_\Delta$ has that property. In 
Proposition~6.1, ii) below, we shall see that equivariant formality 
can be characterized by the freeness of the $A\b$-module 
$IH_{\TT}\b(X)$. Using the 
isomorphism $IH_{\TT}\b(X_\Delta) \cong \E\b(\Delta)$ together with 
the fact that minimal extension sheaves exist on non-rational fans as 
well, this allows to introduce the notion of a (``virtually'') 
equivariantly formal fan in the non-rational 
case, thus eventually leading to a notion of ``virtual'' intersection 
cohomology for arbitrary equivariantly formal fans.

So far, we know that complete fans and $n$-dimensional affine fans are 
of this type. In order to study further examples, we first collect 
some properties pertinent to equivariantly formal fans. 
 
\begin{proposition} \label{prop6.1} 
  \hangindent=1.8\parindent\hangafter=1 
{\rm i)} The $A\b$-module $\IH\b_\TT(\Delta)$ is torsion-free if and 
  only if we have $\Delta$ is purely $n$-dimensional. 
\begin{enumerate} 
\item[ii)] The $A\b$-module $\IH\b_\TT(\Delta)$ is free if and only 
  if $\Delta$ is equivariantly formal. 
\item[iii)] An equivariantly formal fan $\Delta$ is purely $n$-dimensional. 
\item[iv)] If $\Delta$ has an equivariantly formal subdivision 
$\hat\Delta$, then $\Delta$ itself is equivariantly formal. 
\end{enumerate} 
\end{proposition} 
 
We note explicitly that, as a consequence of ii), a notion of 
``virtual'' equivariant formality may be defined even for not 
necessarily rational fans in $\RR^n$ via minimal extension sheaves 
-- no toric varieties are needed. 
 
\begin{proof} i) ``$\Longleftarrow$'': Since $\IH\b_\TT$ is a sheaf, 
we see that we have a natural inclusion 
$$ 
  \IH\b_\TT (\Delta) \;\subset\;\bigoplus_{\sigma 
  \in \Delta^{\max}}\IH\b_\TT(\sigma) 
$$ 
of $A\b$-modules. By Corollary~4.5, each $\IH\b_\TT(\sigma)$ is a 
free $A\b_\sigma$-module. For $\dim\sigma=n$, we have $A\b \cong 
A\b_\sigma$, so the right hand side is a free $A\b$-module. Moreover, 
every submodule of a torsion-free module is again torsion-free. 
 
\noindent ``$\Longrightarrow$'': If $\sigma \in \Delta$ is a maximal 
cone of dimension $d<n$, let $\Delta' := \Delta \setminus \{\sigma\}$. 
The product of some non-zero polynomial function $h \in A\b_\sigma$ 
vanishing on $\partial\sigma$ (such a function can be obtained as a 
product $h = \prod_{\tau \prec_1 \sigma}\ell_\tau$ of non-zero linear 
functions $\ell_\tau \in A^2_\sigma$ with $\ell_\tau|_\tau = 0$) and 
of a non-zero section in $\IH_\TT\b(\sigma)$ yields a non-zero section 
$f \in \IH\b_\TT(\sigma)$ (recall that $\IH\b_\TT(\sigma)$ is a free 
$A\b_\sigma$-module!) that vanishes on $\partial\sigma$. We extend it 
trivially outside of $\sigma$ and thus get a non-trivial torsion 
element, since it is ``killed'' by every non-zero global linear 
function in~$A^2$ vanishing on~$\sigma$. 
 
\noindent ii): ``$\Longleftarrow$'': This is clear by Lemma~4.1, i). 
 
\noindent ``$\Longrightarrow$'': We only sketch the argument, leaving 
details for future exposition: Consider the intersection cohomology 
version of the {\it Eilenberg-Moore spectral sequence\/} (see, e.g., 
\cite[\S~7.2.1]{MCl}) that computes the intersection cohomology of the 
pull back of a bundle. Here we look at the bundle $X_\TT \to B\TT$ and 
take as map the inclusion of a one point set $\{\pt\}$ into $B\TT$. 
Since $B\TT$ is simply connected, the spectral sequence converges: We 
have $E_2^{p,q} \Longrightarrow IH^{p+q}(X)$ with 
$$ 
  E_2^{p,q} \cong \Tor^{p,q}_{H\b(B\TT)}\big(H\b(\pt),IH\b_\TT(X)\big) 
  \cong \Tor^{p,q}_{A\b}\big(\QQ\b,IH\b_\TT(X)\big)\, . 
$$ 
Here $A\b$, $\QQ\b$, and $IH\b_\TT(X)$, respectively, are considered 
as differential graded algebras resp.\ modules with trivial
differential, and we can compare with the classical Tor functors of 
commutative algebra: If $M\b$, $N\b$ are graded $A\b$-modules, the 
corresponding Tor-modules are again in a natural way graded
$A\b$-modules: 
$$ 
  \Tor_p^{A\b}(M\b,N\b) =\bigoplus_{q=0}^\infty \Tor^q_p(M\b,N\b) 
$$ 
\noindent such that 
$$ 
  \Tor^{p,q}_{A\b}(M\b,N\b) \cong \Tor^q_{-p}(M\b,N\b)\, . 
$$ 
Since $IH\b_\TT(X)$ is a free $A\b$-module, we obtain that 
$$
  E_2^{0,\bullet} \cong \QQ\b \otimes_{A\b} IH\b_\TT(X) \cong 
  \overline{IH\b_\TT(X)}
$$ 
and $E_2^{p,q}=\{0\}$ for $p \not= 0$. So in particular, we have  
$E_2^{p,q}=\{0\}$ for $p+q$ odd, and hence also $IH^{p+q}(X)=\{0\}$ 
in that case. 
 
\noindent iii) This follows immediately from i) and ii), since free 
modules are torsion-free. 
 
\noindent iv) According to the Decomposition Theorem of Beilinson,
Bernstein, Deligne, and Gabber, we know that $IH\b(X_\Delta)$ is a 
direct summand of $IH\b(X_{\hat \Delta})$, thus inheriting the 
property that the ``usual'' intersection cohomology vanishes in odd 
degrees. 
\end{proof} 
 
The following example shows that the condition $\Delta^{\max} = 
\Delta^n$ of i) is not sufficient for equivariant formality: 
 
\begin{example} \label{exam6.2} Let $\Delta$ be a fan consisting of two 
equivariantly formal subfans $\Delta_1$ and $\Delta_2$ intersecting 
in a single cone $\tau$. If the codimension of $\tau$ is at least~$2$, 
then $\Delta$ is not equivariantly formal. 
\end{example} 
 
\begin{proof} We intend to prove that the intersection cohomology 
Betti number $Ib_{3}$ is non-zero. With $X := X_{\Delta}$ and $X_i := 
X_{\Delta_i}$, we consider the following part of the exact 
Mayer-Vietoris sequence: 
$$ 
  0 \!\to\! IH^1(X_{\tau}) \!\to\! IH^2(X) \!\to\! IH^2(X_{1}) \oplus
  IH^2(X_{2}) \!\to\! IH^2(X_{\tau}) \!\to\! IH^3(X) \!\to\! 0 \;. 
$$ 
The zeroes at both ends are due to the fact that the toric varieties 
$X_{1}$ and $X_{2}$ are equivariantly formal. 
 
The ``affine orbit splitting'' (0.3) provides an isomorphism $X_\tau 
\cong (\CC^*)^k \times Z_{\tau}$ with $x := \hbox{codim}\,\tau$, 
where $Z_{\tau}$, as a contractible affine $\TT_{\tau}$-toric variety, 
is known to be equivariantly formal. By the K\"unneth formula, we have 
$IH\b(X_\tau) \cong H\b\big((\CC^*)^k\big) \otimes_\QQ IH\b(Z_{\tau})$; 
in particular, we get $Ib^1(X_\tau) = k$ and $Ib^2(X_\tau) = \binom{k}{2} 
+ Ib^2(Z_{\tau})$. By the results of \cite[\S4]{BBFK}, the Betti number 
$Ib^2(X_{\Phi})$ of an arbitrary $d$-dimensional toric variety given by 
a non-degenerate fan~$\Lambda$ (i.e., such that $\Lambda$ spans $V$) is 
determined by the number $a := \#\Lambda^{(1)}$ of rays, namely, we 
have $Ib^2(X_{\Lambda}) = a-d$. Denoting with $a_{i} := 
\#\Delta_{i}^{(1)}$ and $a_{3} := \#\tau^{(1)}$ the respective number 
of rays, we clearly have $\#\Delta^{(1)} = a_{1} + a_{2} - a_{3}$. 
We thus obtain $Ib^2(X_{i}) = a_{i}-n$, $Ib^2(X) = a_{1}+a_{2}-a_{3}-n$, 
and $Ib^2(Z_{\tau}) = a_{3}-(n-k)$. Since the Euler characteristic of 
an exact sequence vanishes, we obtain $Ib^3(X) = k(k-1)/2 > 0 $. 
\end{proof} 
 
We now state a necessary ``topological'' condition for a fan to be 
equivariantly formal, thus obviously providing many examples of toric 
varieties not having that property. 
 
\begin{proposition} \label{prop6.3} If $\Delta$ is an equivariantly 
formal fan, then the complement $N_{\RR} \setminus |\Delta|$ of the 
support $|\Delta|$ is connected. 
\end{proposition} 
 
This is a consequence of the following inequality: 
 
\begin{lemma} \label{lem6.4} The intersection cohomology Betti number 
$Ib^{2n-1}(X_{\Delta})$ satisfies the inequality 
$$ 
Ib^{2n-1}(X_{\Delta}) \ge b_0(N_{\RR} \setminus |\Delta|)-1 \;. 
$$ 
\end{lemma} 
 
\begin{proof} It is known (see, e.g., \cite[3.5]{KaFi}) that the
  natural 
homomorphism 
$$ 
  IH^{2n-1}(X_{\Delta}) \longto H_1^\cld(X_{\Delta}) 
$$ 
is surjective (and even an isomorphism, though we do not need this 
stronger result). To investigate the target, let $\tilde\Delta$ be a 
completion of the fan. The set of cones $\tilde\Delta \setminus 
\Delta$ defines a closed invariant subvariety $\tilde{A}$ of the 
compact toric variety $\tilde{X} := X_{\tilde\Delta}$ that has the 
same number $b_0(\tilde{A})$ of connected components as $N_{\RR} 
\setminus |\Delta|$. Combining that with the isomorphism 
$H_1^\cld(X_{\Delta}) \cong H_1(\tilde{X}, \tilde{A})$ and with the 
exact sequence 
$$ 
  \dots \longto H_1(\tilde{X}, \tilde{A}) \longto H_0(\tilde{A}) 
  \longto H_0(\tilde{X}) \longto H_0(\tilde{X}, \tilde{A}) = 0 
$$ 
eventually yields the following chain of inequalities 
$$ 
  Ib^{2n-1}(X_{\Delta}) \ge b^{2n-1}(X_{\Delta}) \ge 
  b_{1}(X_{\Delta}) = b_{1}(\tilde X, \tilde A) \ge 
  b_{0}(\tilde A) - 1 = b_{0}(N_{\RR} \setminus |\Delta|) 
$$ 
that proves the assertion. 
\end{proof} 
 
Suppose a fan $\Delta$ is obtained from an equivariantly formal fan 
$\Delta_0$ by adding some $n$-dimensional cone together with its 
faces. Then is natural to ask if $\Delta$ is again equivariantly 
formal. The above example~6.2 shows that $|\Delta_0| \cap \sigma$ 
should not be of too small dimension. 
 
Some partial positive results are given by the following propositions. 
 
\begin{proposition} \label{prop6.5} If an equivariantly formal fan 
$\Delta_0$ and an $n$-dimensional cone $\sigma$ intersect in a 
single simplicial cone $\tau \in \Delta$ that is a facet of $\sigma$, 
then the enlarged fan 
$$ 
  \Delta := \Delta_0 \cup \langle \sigma\rangle 
$$ 
is equivariantly formal. 
\end{proposition} 
 
\begin{proof} We have to prove that the toric variety $X := X_\Delta$ 
has vanishing intersection cohomology in each odd degree~$q$ if $X_0 
:= X_{\Delta_0}$ has. To that end, we look at the exact sequence 
$$ 
  \dots \longto IH^{q-1}(X_0) \longto IH^{q}(X, X_0) 
  \longto IH^{q}(X) \longto IH^{q}(X_0) = 0 
$$ 
where the final term vanishes since $q$ is odd. It clearly suffices to 
prove that $IH^{q}(X, X_0)$ vanishes. We may identify this relative 
group with $IH^{q}(X_{\sigma}, X_{\tau})$ by excision. We thus consider 
the analogous exact sequence 
$$ 
  \dots \longto IH^{q-1}(X_{\tau}) \longto IH^{q}(X_{\sigma}, X_{\tau}) 
  \longto IH^{q}(X_{\sigma}) = 0 \,. 
$$ 
As $\tau$ is simplicial, there is an isomorphism $X_{\tau} \cong 
(\CC^{n-1}/F) \times \CC^{*}$, where $F$ is a finite subgroup of 
$(\CC^{*})^{n-1}$ acting diagonally on $\CC^{n-1}$. By the K\"unneth 
formula, we thus have isomorphisms $IH\b(X_{\tau}) \cong 
\QQ\b \otimes_{\QQ} H\b(\CC^{*}) \cong H\b(S^1)$, so $IH^{q}(X_{\tau})$ 
vanishes for each $q > 1$. For $q \ge 3$, these facts immediately yield 
$IH^{q}(X_{\sigma}, X_{\tau}) = 0$. The remaining case $q=1$ follows 
from the fact that the restriction mapping $IH^0(X_{\sigma}) \to 
IH^0(X_{\tau})$ is an isomorphism. \end{proof} 
 
If the cone $\sigma$ is simplicial, we may even allow that it meets 
$\Delta$ in several facets: 
 
\begin{proposition} \label{prop6.6} If an equivariantly formal fan 
$\Delta_0$ and an $n$-dimensional simplicial cone $\sigma$ intersect 
in a subfan $\Lambda$ generated by facets $\tau_{i} \prec_1 \sigma$ 
for $i=1, \dots, \ell$, (i.e., each $\tau_{i}$ is a cone of 
$\Delta_0$), then the enlarged fan 
$$ 
  \Delta := \Delta_0 \cup \langle \sigma \rangle 
$$ 
is equivariantly formal. 
\end{proposition} 
 
\begin{proof} As in the proof of Proposition 6.5, we have to show 
that $IH\b(X,X_0) \cong IH^q(X_{\sigma}, X_{\Lambda})$ vanishes in 
odd degrees. There is an isomorphism 
$$ 
  (X_{\sigma}, X_{\Lambda}) \cong \big(\CC^{n-\ell} \times 
  (\CC^\ell, \CC^\ell \setminus \{ 0\})\big)\bigm/ F \;, 
$$ 
where $F$ is a finite subgroup of $(\CC^{*})^n$ acting diagonally on 
$\CC^n$, such that $\TT \cong (\CC^{*})^n/F$. As passing to the 
quotient by $F$ does not influence the rational (intersection) 
cohomology (see \cite[Thm. II.19.2]{Bre}), we obtain 
$$ 
  IH^q(X_{\sigma}, X_{\Lambda}) \cong H^q(X_{\sigma}, X_{\Lambda}) 
  \cong H^q(\CC^\ell, \CC^\ell \setminus \{ 0 \}) =0 
  \quad\hbox{for}\quad q \not = 0,2\ell \;, 
$$ 
thus proving the assertion.\end{proof} 

Certain fans obtained from complete ones by removing a few 
$n$-dimensional cones are equivariantly formal. This is implied by 
the following result, stated in \cite[Th.~4.2]{Oda$_2$} and 
rephrased here in terms of equivariant formality: 
 
\begin{theorem} \label{ishi} {\rm [Ishida]} The subfan supported 
by the complement of the star of a single ray in a complete fan is
equivariantly formal. 
\end{theorem}

\begin{corollary} \label{cor6.8}  
  \hangindent=1.8\parindent\hangafter=1
{\rm i)} A purely $n$-dimensional fan with convex support is 
equivariantly formal. 
\begin{enumerate} 
\item[ii)] A purely $n$-dimensional fan with open convex complement is  
equivariantly formal. In particular, if $X$ is a complete toric 
variety and $x_\sigma$, a fixed point, then $X \setminus \{x_\sigma\}$ 
is equivariantly formal 
\end{enumerate} 
\end{corollary} 
 
\begin{proof} We choose a vector $v$ such that either $-v$ lies in the 
interior of the support (case~i), or that $v$ lies in the interior of  
the complement (case~ii), and complete the given fan~$\Delta$ by adding  
the new ray $\rho := \QQ_{\ge0}v$ together with all cones of the form  
$\rho + \tau$, where $\tau$ is a cone in the boundary of the support.  
The ``old'' fan then is the complement of the star of the  ``new''  ray. 
\end{proof} 
 
Proposition 6.3 above shows that the completeness hypothesis in the 
second part of~ii) cannot be omitted: If~$\Delta$ is a non-complete fan 
and $\sigma \in  \Delta$, an $n$-dimensional cone, then the support of 
the subfan $\Delta  \setminus \{\sigma\}$ that is defining $X_\Delta 
\setminus \{x_\sigma\}$  has a non-connected complement, since $|\Delta 
\setminus \{\sigma\}|$  is obtained from $|\Delta|$ by removing the 
interior of~$\sigma$, which is separated by the facets from the 
exterior and thus, from the complement of~$|\Delta|$. 
\medskip
 
Using Propostition 6.2 and some homological algebra, we shall prove in 
the companion article \cite{BBFK$_{2}$} the ``topological'' 
characterization of equivariantly formal fans given in the theorem 
below. To formulate it, we use the following notation: For a purely 
$n$-dimensional fan~$\Delta$, we denote with $\partial \Delta$ the 
subfan supported by the topological boundary of~$|\Delta|$ (this is 
the subfan generated by those $(n-1)$-dimensional cones that are 
contained in only one $n$-dimensional cone). Fix a euclidean norm in the 
real vector space $N_\RR$, and denote with~$S$ resp.\ with $\partial S$ 
the intersection of the unit sphere with the support of the fan 
$\Delta$ resp.\ of $\partial \Delta$. Furthermore let $Z := S \cup 
c(\partial S)$ be the compact topological space obtained by patching 
together $S$ and the real cone $c(\partial S)$ over $\partial S$ along 
their respective boundaries.

\begin{theorem} \label{kach} A purely $n$-dimensional fan~$\Delta$ is 
equivariantly formal if and only if the following conditions hold: 
The pair $(S, \partial S)$ has the
real homology of an $(n-1)$-ball modulo its boundary, and $Z$ is a real
homology manifold outside the vertex of the cone $c(\partial S)$.
The last condition is satisfied e.g. if $S$ is a manifold with
boundary. 
\end{theorem}

\typeout{References} 

\typeout{End of document} 
\end{document}